\documentclass{conm-p-l}
\usepackage{amsfonts,amssymb,amscd,amsmath,enumerate,verbatim}

\newcommand{\lra}{\longrightarrow}
\newcommand{\Lra}{\Longrightarrow}
\newcommand{\xra}{\xrightarrow}
\newcommand{\hra}{\hookrightarrow}
\newcommand{\wh}{\widehat}
\newcommand{\up}[1]{{{}^{#1}\!}}

\def\Tor{\operatorname{Tor}} 
\def\Ext{\operatorname{Ext}}
\def\Hom{\operatorname{Hom}}
\def\Ker{\operatorname{Ker}}

\def\Im{\operatorname{Im}}
\def\Soc{\operatorname{Soc}}
\def\Spec{\operatorname{Spec}}
\def\Supp{\operatorname{Supp}}

\def\ch{\operatorname{ch}}
\def\divcl{\operatorname{div}} 
\def\dim{\operatorname{dim}} 
\def\codim{\operatorname{codim}}
\def\embdim{\operatorname{edim}}
\def\depth{\operatorname{depth}}
\def\grade{\operatorname{grade}}
\def\cx{\bullet}
\def\s{\operatorname{s}}
\def\projdim{\operatorname{pd}}
\def\flatdim{\operatorname{fd}}
\def\col{\operatorname{col}}

\def\fpd{\operatorname{fpd}}

\def\crs{\operatorname{crs}}
\def\drs{\operatorname{drs}}
\def\pd{\operatorname{pd}}
\def\id{\operatorname{id}}
\def\cxy{\operatorname{cx}}
\def\curv{\operatorname{curv}}
\def\limsup{\operatorname{limsup}}

\def\fm{{\mathfrak m}}
\def\fn{{\mathfrak n}}
\def\fp{{\mathfrak p}}
\def\fq{{\mathfrak q}}
\def\fel{{l}}

\def\bst{{\boldsymbol t}}
\def\bsx{{\boldsymbol x}}

\def\BQ{{\mathbb Q}}

\theoremstyle{plain}

\theoremstyle{plain}

\theoremstyle{plain}
\newtheorem{thm}{Theorem}[section]

\theoremstyle{definition}
\newtheorem{conj}[thm]{Conjecture}
\newtheorem{ques}[thm]{Question}
\newtheorem{exa}[thm]{Example}
\newtheorem{rem}[thm]{Remark}

\theoremstyle{plain}
\newtheorem{ssthm}{Theorem}[subsection]
\newtheorem{ssprop}[ssthm]{Proposition}
\newtheorem{sscor}[ssthm]{Corollary}
\newtheorem*{kunz}{Kunz's Theorem}
\newtheorem*{cst}{Cohen's Structure Theorem}
\newtheorem*{lcf}{Local Criterion for Finite Flat Dimension}
\newtheorem{fg}[ssthm]{Finite generation}
\newtheorem{fl}[ssthm]{Finite length}

\newtheorem{cxycomp}[ssthm]{Comparison of complexities and curvatures}
\newtheorem{comp}[ssthm]{Comparison of $A$-structures}
\newtheorem{flcomp}[ssthm]{Comparison of lengths}
\newtheorem{AvsB-structure}[ssthm]{Comparison of $A$- and $B$-structures}

\newtheorem{cx/curv}[ssthm]{Complexity and curvature}

\theoremstyle{definition}
\newtheorem{ssrem}[ssthm]{Remark}
\newtheorem{ssdefn}[ssthm]{Definition}
\newtheorem{ssexa}[ssthm]{Example}
\newtheorem{factor}[ssthm]{Factorization}
\newtheorem{compl}[ssthm]{Completion}
\newtheorem{filtr}[ssthm]{Filtration}
\newtheorem*{ssrem*}{Remark}
\newtheorem*{ssdefn*}{Definition}
\newtheorem*{ssntn*}{Notation}
\newtheorem*{proofcite}{{\it{Proof}}}

\begin{document}

\title[Frobenius endomorphism]{The Frobenius endomorphism and \\
homological dimensions}
\author[C.~Miller]{Claudia Miller}
\address{Department of Mathematics, Syracuse University,
Syracuse, NY 13244, U.S.A.}
\email{clamille@msri.org}
\thanks{The author was partly supported by a grant from the 
National Science Foundation.\endgraf 
The article was written during the semesters spent at the 
Universities of Toronto and Nebraska; the author thanks these 
institutions for their hospitality.}
 \maketitle

\section*{Introduction}

In 1969 Kunz \cite{K} proved a fundamental result, connecting the
regularity of a local ring of positive characteristic with the
flatness of its Frobenius endomorphism $\varphi$. This was a first
indication of the important role that $\varphi$ would play in
homological commutative algebra, especially in reflecting  
basic homological properties of the ring.

Some results in Peskine and Szpiro's groundbreaking thesis, announced
in \cite{PS1} in 1969 and published with full proofs in \cite{PS2} in
1973, carry this idea further. Peskine and Szpiro discovered that base
change along $\varphi$ preserves resolutions of modules of finite
projective dimension and then applied this result as a main ingredient
in their proof of the Intersection Theorem. Thus they introduced the
use of the Frobenius map as a tool in solving homological problems in
positive characteristic.  This approach has been developed and used to
prove some major theorems in homological commutative algebra in the
last quarter century.

Herzog in 1974, and others more recently, established various 
converses of Peskine and Szpiro's theorem. Such results extend another
aspect of Kunz's theorem, namely the fact that the Frobenius
endomorphism detects finite homological dimensions. As yet, the full
extent of its power to do so is not completely understood. In this
article we survey what is known about this property of the Frobenius
endomorphism and discuss some questions that remain.

Another major theme of this survey is the behavior of numerical
functions defined by the Frobenius endomorphism that are closely tied
to its homological properties. The study of such functions, once
again, originated in Kunz's Theorem.

In this survey we give, whenever possible, proofs or motivations of
the results we discuss, we record the most precise conclusions that can
be drawn from the proofs, and we recast a few arguments into a
shorter or more illuminating form by using more recent techniques,
sometimes obtaining a new result.

The outline of this paper is as follows: We begin in Section
\ref{kunzsection} with a discussion of the theorem of Kunz mentioned
above. Section \ref{homologicalgens} is devoted generalizations of the
equivalence of regularity with the homological conditions in Kunz's
theorem to the setting of modules.  Section \ref{rodiciosection}
concerns what properties of the ring other than regularity are
reflected by the homological properties of the Frobenius endomorphism.
In Section \ref{numericalgens} we discuss generalizations of the
equivalence with the numerical conditions in the theorem.

Since the situation for complete intersection rings is special, and
indeed much more is known than in general, it is treated separately in
Section \ref{cisection}.

We consider asymptotic versions of the numerical
conditions in Section \ref{asymptotic} and finish with some questions
in Section \ref{questions}. In an appendix we discuss in detail the
properties of the two possible actions of $R$ on the base change of an
$R$-module via the Frobenius endomorphism, and on its derived
functors. While well-known to experts, this material may be useful to
beginners. 

In a different direction, the homological properties of the Frobenius
endomorphism have been used in multiple ways in homological algebra,
not only in proving important theorems, but also in spawning new
constructions that have provided illumination for further steps. For
example, Peskine and Szpiro, Hochster, and P.\ ~Roberts have applied
these to solve central questions in Intersection Theory, cf.\
\cite{BH} and \cite{RobBook} for further reading. The Frobenius
endomorphism has also inspired notions useful in the study of
singularities in geometry, such as $F$-regularity and $F$-rationality,
and has led to the concept of tight closure introduced by Hochster and
Huneke, where many homological properties of the Frobenius
endomorphism are used, cf.\ \cite{Hu} for a systematic exposition.

\subsection*{Conventions}

Before beginning, we summarize some relevant conventions and
terminology. Throughout this paper we let $R$ denote a commutative
Noetherian ring of dimension $d$ and characteristic $p>0$, unless
otherwise indicated. In particular, a local ring is assumed to be
Noetherian as well. When $R$ is local, $\fm$ will always denote its
maximal ideal and $k$ the residue field $R/\fm$.  The embedding
dimension of $R$, denoted $\embdim R$, is the number
$\dim_k\fm/\fm^2$, and the codimension, denoted $\codim R$, is the
number $\embdim R-\dim R$.  We use $\ell_R(-)$, $\projdim_R(-)$ and
$\flatdim_R(-)$ to denote length, projective dimension and flat
dimension, respectively.

The Frobenius endomorphism $\varphi\colon R\to R$ is defined by
$\varphi(r)=r^p$ for $r\in R$. When necessary, we indicate the
relevant ring $R$ by a subscript, writing $\varphi_R$ instead of
$\varphi$. Each iteration $\varphi^n$ defines on $R$ a new structure
of $R$-module, denoted by $\up{\varphi^n}R$, for which $a\cdot
b=a^{p^n}b$.

\section{Regularity}\label{kunzsection}

The ability of the Frobenius map to detect an essential property of
the ring, namely a singularity, was discovered by Kunz in
1969\footnote{Nagata's review \cite{N} of \cite{K} begins as follows:
``Many characterizations of regularity of a local ring are known. The
article gives one which is of a very different type from those known
and which is good only for the case of prime characteristic $p$.''}.
Set $\fm^{[p^n]}=\varphi^n(\fm)R$: this is the ideal generated by the
$p^n$-th powers of any set of generators of $\fm$.

\begin{kunz}[Kunz, {\cite[Thms.~2.1 and 3.3]{K}}]
\label{Kunz}
The following conditions are equivalent for a local ring $R$ 
of characteristic $p$ and dimension $d$.
\begin{enumerate}
\item[\rm{(a)}]
$R$ is regular;
\item[\rm{(b)}]
$\varphi$ is flat;
\item[\rm{(b$^\prime$)}]
$\varphi^n$ is flat for some $n>0$;
\item[\rm{(c)}]
$\ell(R/\fm^{[p]}) = p^{d}$;
\item[\rm{(c$^\prime$)}]
$\ell(R/\fm^{[p^n]}) = p^{nd}$ for some $n>0$.
\end{enumerate}
\end{kunz} 

Kunz's proof proceeds through some variants of the conditions above.
We list these explicitly in order to make specific comparisons
with later results. 

\begin{kunz}[continued]
The following conditions are equivalent to the conditions listed above. 
\begin{enumerate}
\item[\rm{(b$^{\prime\prime}$)}]
$\varphi^n$ is flat for infinitely many $n>0$;
\item[\rm{(b$^{\prime\prime\prime}$)}]
$\varphi^n$ is flat for all $n>0$;
\item[\rm{(c$^{\prime\prime}$)}]
$\ell(R/\fm^{[p^n]}) = p^{nd}$ for infinitely many $n>0$;
\item[\rm{(c$^{\prime\prime\prime}$)}]
$\ell(R/\fm^{[p^n]}) = p^{nd}$ for all $n>0$;
\item[\rm{(d$^{\prime\prime}$)}]
$\ell(R/\fm^{[p^n]}) = p^{ne}$ for infinitely many $n>0$, where $e=\embdim R$.
\end{enumerate}
\end{kunz}

Note that the equivalence of (a) and the various forms of (b) 
does not, in fact, require $R$ to be local since regularity and 
flatness can each be checked locally (cf.\ \cite[Thms.~7.1, 19.3]{Ma}). 
The infinitely many values of $n$ in condition (d$^{\prime\prime}$)
are indeed necessary.

We give an indication of how the proof of the theorem goes; the details
may be found in \cite{K}. Used several times in the proof is the following 
case of Cohen's Structure Theorem. 

\begin{cst}\label{cst}
If $R$ is a complete local ring containing a field, there is a 
ring of formal power series $Q=k[[\bst]]$ on indeterminates $\bst=t_1,
\dots,t_e$, $e=\embdim R$, such that $R=Q/I$ for some ideal $I$ of $Q$.  
If $R$ is regular, $I=0$. 
\end{cst}

The Local Criterion for Flatness is used in Kunz's original proof. We state 
an extension that will be useful in later arguments as well (cf.\
{\cite[Prop.~2.57]{An}}).

\begin{lcf}
Let $R\to S$ be a local homomorphism of local rings 
(of any characteristic), and let $N$ be a finitely generated $S$-module. 
If\, $\Tor_i^R(k,N)=0$, then $\flatdim N<i$.
\end{lcf}

\begin{proof}[Sketch of Proof of Kunz's Theorem]
A first round of implications is\linebreak summarized by the following diagram.
\[
\text{(a)} \Lra \text{(b)} \Lra \text{(b$^{\prime\prime\prime}$)} \Lra 
\text{(b$^{\prime}$)} \Lra \text{(b$^{\prime\prime}$)} \Lra 
\text{(d$^{\prime\prime}$)} \Lra \text{(a)}
\]

For the first implication, let $\iota\colon R\to\wh R$ be the
canonical map to the $\fm$-adic completion $\wh R$, which is
isomorphic to a ring of formal power series $k[[\bst]]$. The map
$\varphi_{\wh R}$ can be factored as $k[[\bst]]\cong k^p[[\bst^p]]\hra
k[[\bst^p]]\hra k[[\bst]]$, where $k^p=\{ x^p\mid x\in k \}$. Note
that $k[[\bst]]$ is free over $k[[\bst^p]]$, and the first map is flat
by the local criterion of flatness:
$\Tor^{k^p[[\bst^p]]}_i(k^p,k[[\bst^p]])$ can be seen to be zero for
$i>0$ by computing it from a Koszul resolution of $k^p$ over
$k^p[[\bst^p]]$.  So, $\varphi_{\wh R}\iota=\iota\varphi_R$ is flat
and then $\varphi_R$ is flat by the faithful flatness of $\iota$.

Of the next three implications, the middle one is trivial and the
others follow from the fact that compositions of flat maps are flat.
To prove that (b$^{\prime\prime}$) implies (d$^{\prime\prime}$), Kunz
shows that if $\varphi^n$ is flat, then the set of $p^n$-th powers of
a set of minimal generators of $\fm$ is independent in the sense of
Lech \cite{L}; properties of independent sets established in \cite{L}
give the desired implication.

Next, using the fact that lengths are preserved under completion and
taking a presentation $k[[\bst]]/I$ for $\wh R$ from Cohen's Structure
Theorem, Kunz observes that condition (d$^{\prime\prime}$) implies
that the lengths of $k[[\bst]]/(I+(\bst^{p^n}))$ and
$k[[\bst]]/(\bst^{p^n})$ agree for infinitely many $n>0$, which is
only possible if $I\subseteq\bigcap_n(\bst^{p^n})=0$ and so $R$ is
regular.

The remaining equivalences are proved in the following order.
\[
\text{(a)} \Lra \text{(c$^{\prime\prime\prime}$)} \Lra \text{(c)} \Lra 
\text{(c$^{\prime}$)} \Lra \text{(c$^{\prime\prime}$)} \Lra 
\text{(b$^\prime$)}
\]

Since lengths are preserved under completion, the first implication
again follows easily from Cohen's Structure Theorem. The next two are
trivial.  That (c$^{\prime}$) implies (c$^{\prime\prime}$) is proved
by showing that an inequality 
\[
\ell(R/\fm^{[p^n]}) \geq p^{nd}
\]
holds for all $R$. 
For the last step, (c$^{\prime\prime}$) implies (b$^\prime$), first
the inequality above implies that $R$ is a domain, and then an
application of Noether's Normalization Theorem and counting of ranks
show that $\up{\varphi^n}R$ is free over $R$.
\end{proof}

Since the regularity of $R$ is equivalent to the condition that the
residue field $k$ has finite projective (equivalently, flat)
dimension, the equivalence of (a) and (b$^\prime$) relates the
homological algebra of the modules $k$ and $\up{\varphi^n}R$. We will
see this theme repeated in other results, most notably in Sections
\ref{homologicalgens}, \ref{rodiciosection} and \ref{cisection}.

The numerical function implicit in condition (c$^{\prime}$) bears
similarities with the Hilbert-Samuel function $\ell(R/\fm^n)$, which
is a polynomial of degree $d$ for $n\gg 0$ with leading coefficient
$e(\fm,R)$ defining the multiplicity of $R$. This has led to the
developments described in Section \ref{asymptotic}.

\section{Finite projective dimension}
\label{homologicalgens}

The flatness condition in parts (b) of Kunz's Theorem can be recast in
the form: $\Tor_i^R(M,\up{\varphi^n}R)=0$ for all finite $R$-modules
and all $i>0$. The equivalence of this condition with (a), the
regularity of $R$, conjures up analogues with the homological behavior
of the residue field $k$. Thinking along these lines has taken two
different routes. 

On the one hand, for a module over any local ring $R$ it is known that
the finiteness of the projective dimension of $M$ is equivalent to the
vanishing of $\Tor_i^R(M,k)$ for one, or for infinitely many, values
of $i$. Thus $k$ always is a test module for finite projective
dimension. In both parts of this section we focus on the extent to
which $\up{\varphi^n}R$ has, or is known to have, similar properties.
The case of complete intersection rings is treated separately in
Section \ref{rigidity} as it requires different machinery.

On the other hand, restrictions on other homological invariants of the
$R$-modules $\up{\varphi^n}R$, such as various homological dimensions,
or certain measures of the size of its resolution, are known to
characterize, just as for $k$, properties of the ring other than
regularity. This is the point of view taken in Section
\ref{rodiciosection}.

\subsection{Finite projective dimension and the Frobenius endomorphism}
\label{PS/Hsection}

The homological manifestations of the Frobenius were further realized
with the result of Peskine and Szpiro in 1969 that base change via the
Frobenius leaves finite free resolutions acyclic. This was their
main tool in proving the Intersection Theorem (\cite[Thm.~1]{PS1},
\cite[Thm.~2.1]{PS2}).

\begin{ssdefn}[Peskine-Szpiro, {\cite{PS1},\cite{PS2}}]
The {\em Frobenius functor\/} $F$, or $F_R$, from the category of 
$R$-modules to itself is given by base change along the 
Frobenius endomorphism $\varphi\colon R\to R$:
\[ 
F(M) = M\otimes_R\up{\varphi}R 
\] 
with the $R$-module structure given through the righthand variable, 
that is, 
\[
r\cdot(m\otimes s)=m\otimes(sr) \quad 
{\text{ for all }}\,r\in R,\, m\in M,\, s\in R.
\]  
\end{ssdefn}

\begin{ssrem}
\label{localisation}
If $M$ is finitely generated, say by $m_1, m_2, \dots, m_r$, then
$F(M)$ is finitely generated, by $m_1\otimes_R 1, 
m_2\otimes_R 1, \dots, m_r\otimes_R 1$. 
For any prime ideal $\fp$ in $R$, since $\varphi^{-1}(\fp)=\fp$, there
are isomorphisms 
\[
(F_R(M))_\fp\cong 
M\otimes_R\up{\varphi}R_\fp \cong 
M_{\varphi^{-1}(\fp)}\otimes_{R_{\varphi^{-1}(\fp)}}\up{\varphi}R_\fp= 
M_\fp\otimes_{R_\fp}\up{\varphi}R_\fp \cong 
F_{R_\fp}(M_\fp).
\] 
It follows that if $M$ is finitely generated, 
$M$ and $F(M)$ have the same support and thus the same dimension.
\end{ssrem}

\begin{ssthm}[Peskine-Szpiro, {\cite[Cor.~2]{PS1}, 
\cite[Thm.~1.7]{PS2}}]
\label{PS} 
Let $R$ be a Noetherian ring of characteristic $p$.
If $M$ is a finitely generated module of finite projective dimension, 
then one has 
\[
\Tor_i^R(M,\up{\varphi}R) = 0
\]
for all $i>0$.
Furthermore, for every prime ideal $\fp$ in $R$ 
\[
\projdim_{R_{\fp}}(F(M))_{\fp}=\projdim_{R_{\fp}}M_{\fp}.
\] 
\end{ssthm}

\begin{proof}
By localizing at a prime minimal in the union of the supports of
$\Tor_i^R(M,\up{\varphi}R)$ for $i>0$, one may assume that for a
minimal free resolution $L_\cx$ of $M$, $F(L_\cx)$ has finite length
homology for $i>0$. Exactness then follows from the Acyclicity Lemma
(\cite[Lemma 1]{PS1}, \cite[Lemma.~1.8]{PS2}), developed for this
purpose, which says that a finite free complex of length at most the
depth of $R$ and with finite length homology modules must be acyclic.
The second assertion follows from the first in view of 
Remarks \ref{localisation} and \ref{Ffree}.
\end{proof}

It follows that if $L_\cx$ is a finite
resolution of $M$ by finitely generated projective modules, then
$F(L_\cx)$ is a projective resolution of $F(M)$. In the local case, 
this translates to the following explicit result by Remark 
\ref{Ffree} below. 

\begin{ssthm}[Peskine-Szpiro, {\cite[Cor.~3]{PS1},
\cite[Thm.~1.13]{PS2}}]
\label{PScomplex} 
Let $R$ be a Noetherian ring of characteristic $p$.
If 
\[
0 \lra L_s \xra{\partial_s} L_{s-1} \xra{\partial_{s-1}} 
\cdots \xra{\partial_1} L_0
\]
is a minimal exact complex of finitely generated free modules, then
the complex 
\[
0 \lra L_s \xra{\partial_s^{[p]}} L_{s-1} \xra{\partial_{s-1}^{[p]}} 
\cdots \xra{\partial_1^{[p]}} L_0
\]
is exact, where $\partial_i^{[p]}$ is given
by the matrix with entries equal to the $p$-th powers of the entries
of the matrix for $\partial_i$ for each $i=1, \dots, s$. 
\qed
\end{ssthm}

The last theorem can be iterated, yielding 
$\Tor_i^R(M,\up{\varphi^n}R) = 0$ for all $i,n>0$ if $M$ has finite
projective dimension.  Iterations of the Frobenius endomorphism play
an important role in what is to come, so we discuss them next.

\begin{ssrem}\label{Ffree}
Let $F^n$, or $F_R^n$, denote $n$ iterations of the Frobenius functor. 
The associativity of tensor products yields an isomorphism 
$F^n(M)\cong M\otimes_R\up{\varphi^n}R$.  Furthermore, 
$F^n(R)\cong R$, and if $f\colon R^s\to R^t$ is a
homomorphism of free $R$-modules given by a matrix $[a_{ij}]$
in some bases, then in the
same bases the homomorphism $F^n(f):R^s \to R^t$ is given by the matrix
$[\varphi^n(a_{ij})]=[a_{ij}^{\ p^n}]$.  Since $F^n$ is right exact, it
follows that if $I$ is an ideal of $R$, then ${F^n_R(R/I)\cong
R/I^{[p^n]}},$ where the {\em bracket power\/} $I^{[p^n]}$ of $I$
denotes the ideal generated by the $p^n$-th powers of any set of
generators of $I$, that is, $I^{[p^n]}=\varphi^n(I)R$.  
\end{ssrem}

\begin{ssrem}
It follows from Theorem \ref{PS} that if $M$ has finite 
projective dimension then so does $F(M)$. 
There are many examples to show that, even if $F(M)$ 
has finite projective dimension, $M$ need not.
For instance, we may take $M$ to be any finitely generated non-free 
module over an Artinian local ring $R$. 
If we choose $n$ so that $\fm^{p^n}=0$, 
then $\fm^{[p^n]}=0$ and $F^n(M)$ is free. 
\end{ssrem}

\begin{ssrem}\label{dp-exs} 
Let $R$ be local. 
By Theorem \ref{PS} and the Auslander-Buchs\-baum equality, whenever $M$
has finite projective dimension, $F^n(M)$ and $M$ have the same depth. 
This may fail otherwise,  
as shown by examples suggested by S.\ Iyengar and A.\ Singh. 

The depth may drop: Let $R=k[[x,y]]/(xy)$ 
where $k$ is a field of characteristic $p$, and let $M=R/I$
where $I=(x)$. Here $M$ has depth one, 
but $F(M)\cong k[[x,y]]/(xy,x^p)$ has depth zero. 
The depth may also rise: Let $R=k[[x,y]]/(x^2)$ and $M=R/I$
where $I=(xy)$. Here $M$ has depth zero, 
but $F(M)\cong k[[x,y]]/(x^2)$ has depth one. 

Combining the ideas in these examples, one may produce an example
where the depth behavior is not monotone: 
Let $R=k[[x,y,z,w]]/(xy,z^{p+1})$ and $M=R/I$
where $I=(x,zw)$. Here $\depth M=1$, 
$\depth F(M)=0$, and $\depth F^2(M)=1$. 
\end{ssrem}

We now discuss a converse of Theorem \ref{PS} that provides a
criterion for finite projective dimension using the Frobenius, and was
the beginning of the realization that the $n$-th Frobenius module
$\up{\varphi^n}R$ may be a test module for finite projective
dimension. We leave out the proof, as a more precise version
by Koh and Lee is proved in Section \ref{kohlee}.

\begin{ssthm}[Herzog, {\cite[Thm.~3.1]{He}}] 
Let $R$ be a Noetherian ring of characteristic $p$, and 
let $M$ be a finitely generated module.
If $\Tor_i^R(M,\up{\varphi^n}R) = 0$ for all $i>0$ 
and infinitely many $n$, then $M$ has finite projective dimension. 
\label{H}
\qed
\end{ssthm} 

To emphasize the analogy with Kunz's Theorem, we summarize these
results. 

\begin{sscor}[Peskine-Szpiro, Herzog]
\label{PS/H}
For a
finitely generated $R$-module $M$ over a
Noetherian ring of characteristic $p$ the following conditions are equivalent.
\begin{enumerate}
\item[\rm{(a)}]
$M$ has finite projective dimension;
\item[\rm{(b$^{\prime\prime}$)}]
$\Tor_i^R(M,\up{\varphi^n}R) = 0$ for all $i>0$ and infinitely many $n$;
\item[\rm{(b$^{\prime\prime\prime}$)}]
$\Tor_i^R(M,\up{\varphi^n}R) = 0$ for all $i,n>0$. 
\qed
\end{enumerate}
\end{sscor} 

This corollary generalizes the equivalence of condition (a) with the 
different versions of condition (b) in Kunz's Theorem. Indeed, by the
Auslander-Buchsbaum-Serre Theorem, condition (a) in Corollary
\ref{PS/H} holds for all finitely generated $R$-modules $M$ if and only if
the ring $R$ is regular. Similarly, condition (b$^{\prime\prime}$)
(resp., (b$^{\prime\prime\prime}$)) in Corollary \ref{PS/H} holds for
all finitely generated $R$-modules $M$ if and only if $\up{\varphi^n}R$ is
flat over $R$ for infinitely many $n$ (resp., for all $n$ ). As
discussed in Proof of Kunz's Theorem, the other implications
between its various conditions (b) are trivial because compositions of
flat maps are flat.

The reader may have noticed that Corollary \ref{PS/H} leaves open the
question of whether the condition on the vanishing of
$\Tor_i^R(M,\up{\varphi^n}R)$ for a single $n$ can be included,
corresponding to condition (b$^{\prime}$) in Kunz's Theorem.  More
specifically, one may ask: For how many values of $i$ and $n$ does one
need $\Tor_i^R(M,\up{\varphi^n}R)$ to vanish to ensure that $M$ has
finite projective dimension? Recently there has been progress on this
question; we describe it in Sections \ref{kohlee} and \ref{rigidity}.

\subsection{Finitistic criteria for finite projective dimension}
\label{kohlee}

To discuss a first criterion, we introduce some invariants defined by
Koh and Lee in 1998, \cite{KL1}. Their ideas grew out of techniques
used by Burch in \cite{B} and Hochster in \cite{Ho3}, and present
implicitly in Herzog's proof of Theorem \ref{H} above. 

In this section $R$ is local and $M$ and $N$ denote finitely
generated $R$-modules.

We will give the most precise statements that can be made from the
proofs in \cite{KL1}. Koh and Lee also develop a dual set of
definitions and results using injective resolutions; we do not
list them here, but mention them briefly at the end of this subsection.

Alongside the definitions given by Koh and Lee, we introduce bracket
power versions of their invariants in characteristic $p$ by replacing
the regular powers $\fm^t$ of $\fm$ with its bracket powers
$\fm^{[p^n]}$.  Except where indicated, their results hold in any
characteristic; for statements involving bracket powers, it is
implicitly assumed that $R$ has characteristic $p$.

\begin{ssdefn}[\cite{KL1}] 
For a homomorphism $\partial\colon R^n\to R^m$, set 
\setlength\arraycolsep{2pt}
\begin{eqnarray*}
\col(\partial) &=& \inf\{t\geq 1\mid  
\pi\circ\partial\,(R^n)\nsubseteq \fm^t 
\text{ for each epimorphism }\pi: R^m\to R\};\\
\col_p(\partial) &=& \inf\{p^n\geq 1\mid  
\pi\circ\partial\,(R^n)\nsubseteq \fm^{[p^n]}
\text{ for each epimorphism }\pi: R^m\to R\}.
\end{eqnarray*}

Thus, $\col(\partial)$ is the least $t$ such that, for 
any\footnote{In {\cite{KL1}} it is not made explicit whether
``some" or ``any" was intended; it was the latter.} $m\times n$
matrix $[a_{ij}]$ representing $\partial$ with respect to some bases
of $R^n$ and $R^m$, each row\footnote{Koh and Lee \cite{KL1} multiply
matrices on the right, so our rows are their columns.}  of $[a_{ij}]$
contains an element outside $\fm^t$; the analogous statement holds for
$\col_p(\partial)$ using bracket powers of $\fm$.

For an $R$-module $M$, let $\partial_i^M$ denote the $i$-th map in a minimal
free resolution of $M$, and set 
\[
\col(M)=
\begin{cases}
\inf\{\col(\partial_i^M)\mid i>1+\depth R\}
&\text{if } \projdim M=\infty;\\
1 &\text{if } \projdim M<\infty.\\
\end{cases}
\]
An invariant $\col_p(M)$ is defined similarly. 
\end{ssdefn}

\begin{ssdefn}[\cite{KL1}]
For an $R$-module $N$, let
$\Soc(N)=(0:\fm)_N$ denote the socle of $N$ and set
\[
\s(N)=\inf\{ t>0\mid  \Soc(N)\nsubseteq\fm^tN \}.
\]
An invariant $\s_p(M)$ is defined similarly.
\end{ssdefn}

\begin{ssdefn}[{\cite{KL1}, \cite{KL2}}]\label{kohleedefn} For every
local ring $(R,\fm)$, set\footnote{The intended association is ``column"
for $\col$, see the preceding footnote, ``finite projective dimension"
for $\fpd$, and ``cyclic by a regular sequence" for $\crs$.}
\setlength\arraycolsep{2pt}
\begin{eqnarray*}
\col(R)&=&
\sup\{\col(M)\mid  M \text{ is an }R\text{-module}\};\\
\fpd(R)&=&\inf\{ \s(N)\mid  N \text{ is an }R\text{-module with }  
\projdim N <\infty \text{ and }\depth N =0 \};\\
\crs(R)&=&\inf\{ \s(R/({\bsx}))\mid 
{\bsx}=x_1, \dots, x_r\text{ is a maximal R-sequence }\}.
\end{eqnarray*}
Define $\col_p(R)$, $\fpd_p(R)$, and $\crs_p(R)$ similarly, by using 
$\col_p(M)$ and $\s_p(N)$ in place of $\col(M)$ and $\s(N)$, respectively.
\end{ssdefn} 

To relate these invariants, we need the next remarkably powerful result 
that has a very simple proof. 

\begin{ssprop}[Koh-Lee, {\cite[Prop.~1.2]{KL1}}]
\label{KLlemma}
If $N$ is an $R$-module of depth $0$ such that $\Tor_i^R(M,N)=0$
for some $i$ with $0< i < \projdim M$, 
then each row of $\partial_{i+1}^M$ 
contains an element outside of $\fm^{\s(N)}$
(resp., an element outside of $\fm^{[\s_p(N)]}$).
\end{ssprop}

\begin{proof} Let $(L_\cx,\partial_\cx)$ be a minimal free resolution
of $M$.  By hypothesis, tensoring $L_\cx$ with $N$ gives a complex
that is exact in the $i$-th spot. If the entire $j$-th row of
$\partial_{i+1}^M$ were contained in $\fm^{\s(N)}$, then the $j$-th
component of every element in $\Im (\partial_{i+1}\otimes_R N)=\Ker
(\partial_i\otimes_R N)$ would be in $\fm^{\s(N)}N$. But since $L_\cx$
is minimal, $L_i\otimes_R\Soc N\subseteq \Ker (\partial_i\otimes_R N)$
and so one would have $\Soc N\subseteq \fm^{\s(N)}$, a contradiction. 
An analogous proof works for the statement involving bracket powers. 
\end{proof}

Applying Proposition \ref{KLlemma} to the modules $N$ used in the
definitions of $\fpd(R)$ and $\fpd_p(R)$ gives the following
corollary; indeed, choose $N$ with $\projdim N <\infty$ and 
$\depth N =0$ so that $\s(N)=\fpd(R)$ (or $\s_p(N)=\fpd_p(R)$). 
Since $\projdim N\leq\depth R$, if $i-1>\depth R$, then 
$\Tor_i^R(M,N)=0$, and we are done by the proposition. 

\begin{sscor}[Koh-Lee, {\cite[Thm.~1.7]{KL1}}] 
\label{KLresns}
If $M$ has infinite projective dimension, then for all $i>1+\depth R$ 
each row of $\partial_i^M$ contains an element outside of
$\fm^{\fpd(R)}$ (resp., an element outside of $\fm^{\fpd_p(R)}$).
\qed
\end{sscor}

This immediately yields a comparison between the invariants in
Definition \ref{kohleedefn}.

\begin{sscor}[Koh-Lee, {\cite[Prop.~1.5]{KL1}}]\label{invariants}
There are inequalities 
\begin{xxalignat}{3}
& &\col(R) &\leq \fpd(R) \leq \crs(R) < \infty;&&\\
&{\hphantom{\square}}& \col_p(R)& \leq \fpd_p(R) \leq \crs_p(R) < \infty.
&&\square
\end{xxalignat}
\end{sscor}

In fact, when the residue field $k$ of $R$ is infinite, the numbers in
the first row are conjectured to be equal (\cite{KL2}), but this can
fail if $k$ is finite (\cite[Ex.~4.8]{KL2}). These numbers are always
equal if $R$ has depth zero, (\cite[Prop.~1.7]{KL2}).

We give examples to show that in general $\crs(R)$ and $\crs_p(R)$ 
are incomparable, and so neither statement in Proposition \ref{KLlemma} 
or in Corollary \ref{KLresns} implies the other.  

\begin{ssexa}
\label{examples} 
Consider the ring $R=k[x,y]/(x^2,y^2)$, where $k$ is a field.  Note
that $R$ is zero-dimensional with $\Soc(R)=(xy)$.  If $k$ has
characteristic 2, then $\Soc R \subseteq \fm^2$ and 
$\Soc R \nsubseteq\fm^3$, but
$\Soc R \nsubseteq \fm^{[2]}$. So, we see that  $\crs(R)=3$ and 
$\crs_p(R)=2$.  However, if $k$ has 
characteristic 5, then $R$ satisfies $\crs(R)=3$ and $\crs_p(R)=5$.  
\end{ssexa}

We now state Koh and Lee's result, which is a significant
strengthening of Theorem \ref{H}, as it gives a finitistic criterion
for finite projective dimension. Koh and Lee state it with the
condition $p^n \geq \col_p(R)$, but this number is not usually
available. The strength of the theorem lies in the fact that, due to
the inequality in Corollary \ref{invariants}, it can be applied by
ensuring that $p^n \geq \crs_p(R)$, or that
$\Soc(R/({\bsx}))\nsubseteq\fm^{[p^n]}/({\bsx})$ for some $R$-sequence
${\bsx}$.  Excepting this remark, the following statement is
the most precise conclusion that one obtains from their proof; it
corrects misprints in \cite{KL1} and \cite{AM}. It is clear that the
bracket invariant yields a stronger (or equal) result, so we give only
this statement.

\begin{ssthm}[Koh-Lee, {\cite[Prop.~2.6]{KL1}}]
\label{KL}
Suppose that $R$ is a local ring of characteristic $p$, and 
let $n$ be an integer such that $p^n \geq \crs_p(R)$. 
Let $M$ be a finitely generated $R$-module. 
If $\Tor_i^R(M,\up{\varphi^n}R) = 0$ for $\depth R +1$ consecutive 
values of $i>0$, then $M$ has finite projective dimension. 
\end{ssthm}

\begin{proof} 
Let $(L_\cx,\partial_\cx)$ be a minimal free resolution of $M$, 
and let $r$ denote the depth of $R$.  
If $\Tor_i^R(M,\up{\varphi^n}R) = 0$ for $\ell+1 \leq i \leq \ell+r+1$, 
then tensoring $L_\cx$ with 
$\up{\varphi^n}R$ and truncating gives a minimal acyclic complex
\[
L_{\ell+r+2} \xra{F^n(\partial_{\ell+r+2})}
L_{\ell+r+1} \lra
\cdots \lra 
L_{\ell+1} \xra{F^n(\partial_{\ell+1})} 
L_{\ell} \lra 
0
\]
of length $r+2$
whose maps have entries in $\fm^{[p^n]}\subseteq\fm^{\crs_p(R)}$.
Because of the inequality in Corollary \ref{invariants},
this contradicts the definition of $\col_p(R)$ unless $L_{\ell+r+1}=0$. 
\end{proof}

For rings of depth zero and $n\geq\log_p(\crs_p(R))$,
this already gives the strongest possible version of Herzog's result.

\begin{sscor}
\label{KLdepth0}
Let $n$ be such that $\Soc(R)\nsubseteq\fm^{[p^n]}$. 
If $\Tor_i^R(M,\up{\varphi^n}R) = 0$ for some $i>0$, 
then $M$ has finite projective dimension.
\qed
\end{sscor}

In particular, for rings $R$ with $\Soc(R)\nsubseteq\fm^{[p]}$, such as 
the rings in Remark \ref{examples}, there is no restriction on $n$. 

\begin{ssrem}
When $R$ is Artinian, the proof of Theorem \ref{KL} trivializes: 
Its hypothesis becomes $\fm^{[p^n]}=0$; hence in the complex 
$L_\cx\otimes_R \up{\varphi^n}R$ all boundary maps are zero, 
and exactness in the $i$-th spot means $L_i=0$.
\end{ssrem}

Now suppose that $R$ is Cohen-Macaulay.  Using modules of finite
injective dimension, Koh and Lee develop a set of definitions and
results in \cite{KL1} dual to those above.  The dual theory can give a
slightly stronger result than the one in Theorem \ref{KL}. We 
sketch their proof, modifying it for bracket powers. 
Again, we give the effective version involving the bracket version of 
the invariant $\drs(R)$.
It can be defined as the ``bracket Loewy length" of $R$, namely, 
\[
\drs_p(R)=
\inf\{ p^n\geq 1\mid \fm^{[p^n]}\subseteq({\bsx}), \text{ for some maximal 
$R$-sequence }{\bsx}\}.
\]
 From this definition it is clear that $\crs_p(R)\leq \drs_p(R)$, 
and it can be shown that equality holds if $R$ is Gorenstein
(by an argument like that for \cite[Prop.\ 4.1]{KL2}). 

\begin{ssthm}[Koh-Lee, {\cite[Prop.~2.6]{KL1}}]
\label{KLdual}
Suppose that $R$ is Cohen-Macaulay and has positive dimension, and 
let $n$ be such that $p^n \geq \drs_p(R)$. 
If\, $\Tor_i^R(M,\up{\varphi^n}R) = 0$ for $\depth R$ consecutive 
values of $i>0$, then $M$ has finite projective dimension. 
\end{ssthm}

\begin{proof}
As in the proof of Theorem \ref{KL}, 
one obtains a minimal acyclic complex 
\[
G \xra{}
L_{\ell+r+1} \xra{F^n(\partial_{\ell+r+1})}
\cdots \lra 
L_{\ell+1} \xra{F^n(\partial_{\ell+1})} 
L_{\ell} \lra 
0
\]
but with $L_{\ell+r+2}$ replaced by another
free $R$-module $G$ to retain exactness at $L_{\ell+r+1}$.
Let $C$ be the cokernel of $F^n(\partial_{\ell+1})$.

Choose a maximal $R$-sequence ${\bsx}$ so that
$\fm^{[p^n]}\subseteq({\bsx})$ as in the definition of $\drs_p(R)$.
The finitely generated $R$-module $T=\Hom_R(R/({\bsx}),E)$, where $E$
is the injective hull of $k$ over $R$, has depth zero and finite
injective dimension.  Since $\Soc(T)$ is one-dimensional, the condition 
$\Soc(T)\nsubseteq\fm^{[p^n]}T$ holds if and only if $\fm^{[p^n]}T=0$. 
Thus $p^n\geq \s(T)$. 

Since $\Ext_R^{r+1}(C,T)=0$, 
applying $\Hom(-,T)$ to the complex above leaves it exact at $L_{\ell+r+1}$. 
Next, an argument exactly as in the proof of Proposition 
\ref{KLlemma} yields that the map to the 
left of $\Hom(L_{\ell+r+1},T)$ has some entries not in $\fm^{[p^n]}$. 
Since the arrows have been reversed, that map is actually the 
$T$-dual of $F^n(\partial_{\ell+r+1})$, and we arrive at a contradiction. 
\end{proof}

\begin{sscor}
\label{KLdepth1}
Suppose that  $R$ is a one-dimensional Cohen-Macaulay ring, and 
let $n$ be such that $p^n \geq \drs_p(R)$. 
If $\Tor_i^R(M,\up{\varphi^n}R) = 0$ for some $i>0$, 
then $M$ has finite projective dimension.
\qed
\end{sscor}

The dual theory also yields a result for injective dimension analogous
to the one given in Theorem \ref{KL}, namely that finiteness of the
injective dimension of a module $M$ can be detected from the vanishing
of $\Ext^i_R(\up{\varphi^n}R,M)$ for $\depth R+1$ consecutive values
of $i>0$ and some $n\geq \crs_p(R)$.  It is interesting to compare
this to an earlier theorem of Goto: if, for some $n>0$,
$\Hom_R(\up{\varphi^n}R,R)\cong\up{\varphi^n}R$ and
$\Ext^i_R(\up{\varphi^n}R,R)=0$ for $1\leq i \leq \depth R$,
then $R$ is Gorenstein \cite{Go}.

\section{Homological dimensions}\label{rodiciosection}

In the previous section, we saw how the Frobenius endomorphism behaves
with respect to modules of finite projective dimension and how it can
be used to detect finite projective dimension. Since regularity of a
ring is equivalent to the finiteness of its global dimension, these
results yield that regularity is equivalent to finite flat dimension
of all iterations of the Frobenius endomorphism. The stronger version
involving only one iteration is discussed in this section.

Recently, the question of what other properties of a ring of
characteristic $p$ can be detected by the finiteness of certain
homological dimensions of the Frobenius endomorphism has been
considered.  In the first subsection we survey results for the
homological dimensions: flat dimension, injective dimension,
CI-dimension, and G-dimension; in the second we discuss a related
asymptotic homological condition. We omit most proofs in this 
section, but point the reader to the appropriate references.

\subsection{Detection of properties of the ring from the structure 
of the Frobenius endomorphism}
\label{homdims}

We begin with the property of regularity: In 1988 Rodicio used
Andr\'e-Quillen homology in an essential way to strengthen a crucial
implication, (b) $\implies$ (a), in Kunz's Theorem. Although not
stated explicitly in \cite{He}, it follows also from the earlier
result of Herzog, namely Theorem \ref{H} above.  Koh and Lee have shown that
it follows easily from their results as well, discussed in Section
\ref{kohlee}; this is the proof that we give here. 

\begin{ssthm}[Rodicio, {\cite[Thm.~2]{Ro}}]
\label{rodicio}
Let $R$ be a local ring of characteristic $p$.
If $\flatdim_R(\up{\varphi^{n}}R)<\infty$ for some $n>0$, 
then $R$ is regular. 
\end{ssthm}

\begin{proofcite}
[{\cite[Prop.~2.6]{KL1}}]
If ${\varphi^n}$ has finite flat dimension, then so do its
self-compositions ${\varphi^{nt}}$, as can be shown using 
the Cartan-Eilenberg spectral sequence. Take $t$ such that
$p^{nt}\geq\crs_p(R)$.  Since $\Tor_i^R(k,\up{\varphi^{nt}}R)=0$ for
all $i>\flatdim_R(\up{\varphi^{nt}}R)$, Theorem \ref{KL} implies that
$k$ has finite projective dimension, so $R$ is regular.  
\qed
\end{proofcite}

By the Local Criterion of Finite Flat Dimension, cf.\
\S\ref{kunzsection}, Theorem \ref{rodicio} can be restated as a
nonvanishing statement for non-regular rings.

\begin{sscor}
\label{non-vanishing}
If $\Tor_i^R(k,\up{\varphi^n}R)= 0$ for some 
$i>0$ and some $n>0$, then $R$ is regular. 
\qed
\end{sscor}

It turns out that regularity is also detected by the injective
dimension of the Frobenius endomorphism.

\begin{ssthm}[Avramov-Iyengar-Miller, {\cite{AIM}}] \label{kunzid}
Let $R$ be a local ring of characteristic $p$. 
The ring $R$ is regular if and only if 
$\id_{R}(\up{\varphi^n}R) < \infty$ for some $n>0$.
\end{ssthm}

While flat or injective dimension are defined directly, other properties
of homomorphisms, such as Gorenstein and the Cohen-Macaulay properties,
as well as the corresponding homological dimensions, are defined via
a Cohen factorization. This enables one to consider a finite (in fact,
surjective) homomorphism and use the definition for finite modules.

\begin{ssdefn}[Avramov-Foxby-Herzog, \cite{AFH}]
Let $\alpha: (A,\fm) \to (B,\fn)$ be a local homomorphism of local rings;
let $\iota\colon B\to \wh B$ be the canonical inclusion of $B$ into its
completion. A {\em Cohen factorization\/} of $\iota\alpha$ 
is a factorization $A\xra{\dot\alpha} A' \xra{\alpha'} \wh B$ such that 
the map $\dot\alpha$ is flat, the ring $A'$ is complete, the ring 
$A'/\fm A'$ is regular, and the map $\alpha'$ is surjective. 
A {\em Cohen factorization\/} always exists, but is not necessarily unique, 
although any such factorization can be reduced to a minimal one, 
for which $\embdim A'/\fm A'= \embdim B/\fm B$ holds. 
\end{ssdefn}

The results below use homological dimensions whose global finiteness
characterizes certain properties of the ring.  The terminology derives
from the fact that a ring $A$ is complete intersection (respectively,
Gorenstein, Cohen-Macaulay) if and only if the CI-dimension (respectively,
G-dimension, CM-dimension) of every finitely generated $A$-module
is finite.

\begin{ssdefn}[Avramov-Gasharov-Peeva, \cite{AGP}]
A finite $A$-module $M$ is said to have {\em finite CI-dimension\/}
if there is a local flat homomorphism $A\to A'$ and a surjective 
homomorphism $A'\gets Q$ with kernel generated by a regular sequence 
such that $\projdim_Q(M\otimes_A A')<\infty$. 
\end{ssdefn}

Blanco and Majadas proved that the Frobenius endomorphism can be used
to detect the complete intersection property (cf.\ \S 5 for the definition ).

\begin{ssthm}[Blanco-Majadas, {\cite[Prop.~1]{BM}}]\label{cidim}
Let $R$ be a local ring of characteristic $p$.  The ring $R$ is complete
intersection if and only if for some $n>0$ and for some (equivalently,
for any) Cohen factorization $R\to R'\to\wh R$ of $\varphi^n$, the module
$\wh R$ has finite CI-dimension over $R'$.  \end{ssthm}

The proof uses techniques similar to those of Rodicio's proof of
Theorem \ref{rodicio}.

\begin{ssdefn}[Auslander-Bridger, \cite{AB}]
A finite $A$-module $M$ is said to have {\em finite G-dimension\/}
if it has a finite resolution by finite modules $G_n$ such that 
each $G_n$ reflexive and satisfies $\Ext_i^R(G_n,R)=\Ext_i^R(G_n^*,R)=0$
for all $i>0$. 
\end{ssdefn}

Iyengar and Sather-Wagstaff have shown that Frobenius endomorphism can
be used to detect the Gorenstein property as well. A proof was given
independently by Takahashi and Yoshino \cite{TY} in the case that
$\varphi$ is finite.

\begin{ssthm}[Iyengar-Sather-Wagstaff, {\cite{ISW}}]\label{gordim}
Let $R$ be a local ring of characteristic $p$.  The ring $R$ is Gorenstein
if and only if for some $n>0$ and for some (equivalently, for any)
Cohen factorization $R\to R'\to\wh R$ of $\varphi^n$, the module $\wh R$
has finite G-dimension over $R'$.
\end{ssthm}

Gerko has introduced a notion of CM-dimension, which combines features
of both definitions above, as follows:

\begin{ssdefn}[Gerko, \cite{Ge}]
A finite $A$-module $M$ is said to have {\em finite CM-dimension\/}
if there is a local flat homomorphism $A\to A'$ and a surjective
homomorphism $A'\gets Q$ such that both $A'$ and $M\otimes_A A'$
have finite G-dimension over $Q$, and $\Ext^i_Q(A',Q)=0$ for $i<\grade_Q(A')$.
\end{ssdefn}

The finiteness of CM dimension of all finite modules  global characterizes
the Cohen-Macaulayness of a ring. When the Frobenius module is finite,
Takahashi and Yoshino have proved that the finiteness of its CM-dimension
suffices.

\begin{ssthm}[Takahashi-Yoshino, {\cite{TY}}]\label{cmdim}
Let $R$ be a local ring of characteristic $p$ and suppose that $k$ is 
perfect.  The ring $R$ is Cohen-Macaulay if and only if for $n\gg 0$
the $R$-module $\up{\varphi^n}R$ has finite CM-dimension.  
\end{ssthm}

\subsection{Detection of properties of the ring from asymptotic data on 
the Frobenius endomorphism}
\label{numerprops}

When $R$ is not regular, so that $\Tor_i^R(k,\up{\varphi^n}R)\neq 0$
for all $i>0$ by Corollary \ref{non-vanishing}, we examine the 
asymptotic behavior of the numerical function 
\[
i \mapsto \ell_{\varphi^n}\big(\Tor_i^R(k,\up{\varphi^n}R)\big),
\]
where lengths over $\varphi^n$ are defined as in 
Appendix \ref{fgsection}. 

The complexity measures the degree of the least polynomial bound of
this function; it is defined as follows, cf.\ Appendix \ref{alphacx}.

\begin{ssdefn}\label{frobcxdef}
For a local ring $R$, the {\em complexity\/} of $R$ over $\varphi^n$, 
denoted $\cxy_{\varphi^n}\!R$, is the least non-negative integer $t$ 
with the property that 
\[
\ell_{\varphi^n}\big(\Tor_i^R(k,\up{\varphi^n}R)\big) \leq \beta\, i^{t-1}
\]
for some $\beta\in {\mathbb{R}}$ and all $i\gg0$; if no such $t$ exists, 
then $\cxy_{\varphi^n}\!R=\infty$.
\end{ssdefn}

The curvature measures the exponential rate of growth of the numerical
function above; it is defined as follows, cf.\ Appendix \ref{alphacx}.

\begin{ssdefn}\label{frobcurvdef}
For a local ring $R$, the {\em curvature\/} of $R$ 
over $\varphi^n$ is defined as 
\[
\curv_{\varphi^n}\!R=
\limsup_i \sqrt[i]{\ell_{\varphi^n}\big(\Tor_i^R(k,\up{\varphi^n}R)\big)}
\]
\end{ssdefn}

The main result of \cite{AIM} is that this function has maximal growth
(cf.\ the inequalities (\ref{kextremal}) in Appendix \ref{alphacx}): 

\begin{ssthm}[Avramov-Iyengar-Miller, {\cite{AIM}}] \label{extremal}
Let $R$ be a local ring of characteristic $p$. For any $n>0$, 
there are equalities
\[
\cxy_{\varphi^n}\!R = \cxy_R k
\qquad\text{and}\qquad
\curv_{\varphi^n}\!R =\curv_R k. 
\]
\end{ssthm}

In particular, this yields a characterization of the complete
intersection property in terms of a purely numerical condition, which
is weaker than the structural condition given in Theorem \ref{cidim}.
See \cite{AIM} for details on how these two results are related, the
main link being that the complexity and curvature can also be computed
via a Cohen factorization of $\varphi^n$. 

\begin{ssthm}[Avramov-Iyengar-Miller, {\cite{AIM}}] \label{cichar}
Let $R$ be a local ring of characteristic $p$. 
The following conditions are equivalent. 
\begin{enumerate}
\item[\rm{(a)}]
$R$ is complete intersection.
\item[\rm{(b$'$)}]
$\cxy_{\varphi^n}\!R < \infty$ for some $n>0$.
\item[\rm{(b$''$)}]
$\cxy_{\varphi^n}\!R = \codim R$ for all $n>0$.
\item[\rm{(e$'$)}]
$\curv_{\varphi^n}\!R \leq 1$  for some $n>0$. 
\item[\rm{(e$''$)}]
$\curv_{\varphi^n}\!R \leq 1$  for all $n>0$. 
\end{enumerate}
\end{ssthm}

Theorem \ref{cichar} follows directly from Theorem \ref{extremal} in
view of an analogous set of characterizations of complete intersection
rings in terms of the complexity and curvature of the residue field.

Since $R$ is regular if and only if $R$ is complete intersection of
codimension zero, Theorem \ref{cichar} extends that part Kunz's Theorem
from Section \ref{kunzsection}, which states that the regularity of
a local ring $R$ is equivalent to the flatness of $\varphi^n$ for any
$n>0$, and also Rodicio's Theorem \ref{rodicio}.

\begin{ssrem}
The results in Section \ref{rodiciosection} again reflect the 
strong similarities of the module $\up{\varphi^n}R$ and the 
residue field $k$ in their roles as test modules for certain homological 
properties of the ring $R$. Each of Theorems \ref{rodicio}, 
\ref{cidim}, \ref{gordim}, \ref{cmdim} and \ref{cichar} are standard results 
for the module $k$ in place of the module $\up{\varphi^n}R$:  
the finiteness of the flat dimension, CI-dimension,
G-dimension, CM-dimension, or complexity, respectively, of $k$  
implies that $R$ is regular, complete intersection, Gorenstein, 
Cohen-Macaulay, or complete intersection, respectively. 
\end{ssrem}

\section{A numerical condition}
\label{numericalgens}

In this section we consider generalizations of Kunz's equivalence of
the regularity of a local ring $R$ with the numerical condition:
\begin{enumerate}
\item[\rm{(c$^{\prime\prime\prime}$)}]
$\ell(R/\fm^{[p^n]}) = p^{nd}$ for all $n>0$. 
\end{enumerate}
Since $R/\fm^{[p^n]} \cong F^n(k)$, this condition can
be expressed in terms of Frobenius functors in the form 
$\ell(F^n(k)) = \ell(k) p^{nd}$.

In this section, $M$ denotes a module of finite length; 
as explained in the Appendix, the length of $F^n(M)$ is then finite 
for all $n\geq 0$.

\begin{rem}\label{lengthsregular}
When $R$ is regular the Frobenius functors are exact, so an easy induction on
the length of $M$ shows that (c$^{\prime\prime\prime}$) is equivalent
to the following seemingly stronger property:
$\ell(F^n(M)) = \ell(M) p^{nd}$ for all $R$-modules $M$ of finite length
and all $n>0$. 
\end{rem} 

In 1974, Peskine and Szpiro gave the first result generalizing Kunz's
numerical characterization of regularity.  

\begin{thm}[Peskine-Szpiro, {\cite[Thm.~2]{PS3}, \cite[\S 3]{Sz}}]
\label{PSgraded}
Let $R=\bigoplus_{i\geq 0}R_i$ be a graded Noetherian ring  
generated by $R_1$ such that $R_0$ is Artinian, and let $M$ be 
a graded $R$-module of finite length. 
If $M$ has finite projective dimension, then 
$\ell(F^n(M)) = \ell(M) p^{nd}$ for all $n\geq 0$.  
\end{thm}

Based on this result, Szpiro \cite{Sz} conjectured that the equality
above holds for all modules of finite length over a local ring (and
predicted a version of the Chern theory discussed in Section
\ref{chern}). In 1985, Dutta showed this is indeed the case for
several classes of rings, namely, complete intersection rings, cf.\
Theorem \ref{Duttaci} (ii), and rings of small dimension.

\begin{thm}[Dutta, {\cite[1.14]{D1}}]
\label{Duttacmgor}
Let $R$ be a Cohen-Macaulay local ring of dimension at most two
or a Gorenstein local ring of dimension three, and 
let $M$ be an $R$-module of finite length. 
If $M$ has finite projective dimension, 
then $\ell(F^n(M)) = \ell(M) p^{nd}$ for all $n\geq 0$.
\end{thm}

For Cohen-Macaulay rings of dimension at most two the result is
deduced from the fact that the Vanishing Conjecture holds over such
rings (cf.\ \cite[1.13]{D1}).  Another proof was provided by the
development of the theory of local Chern characters, cf.\ Section
\ref{chern} for more details.

In the case of complete intersection rings much more is known. 
This is the topic of Section \ref{cinumerical}. In fact, we will 
see that over complete intersection rings the strongest possible
generalization of Kunz's numerical criterion holds, namely that given
by Theorems \ref{Duttaci} and \ref{length}.  However, in general it
fails completely. 

\begin{exa}\label{exs} 
The first example with 
\[
\ell(F^n(M))\neq \ell(M)p^{nd}
\]
for a module of finite projective dimension was given by P.\ Roberts
in 1989, \cite[\S 4]{RobMSRI}. He used a pair of modules $(M',N')$,
constructed  by Dutta, Hochster and McLaughlin \cite{DHM} with a negative
Serre intersection multiplicity over a three-dimensional hypersurface
$R^\prime$ such that $R'$-module $M'$ has finite length and finite
projective dimension.  First, Roberts constructed a finite extension ring
$R$ of $R'$ carefully chosen to exploit the pathological behavior of the
modules (in particular, the class of $R$ is equal to the class of $N'$
in the reduced Grothendieck group of $R^\prime$).  After a base change
of $M'$ to $R$, which is a three-dimensional Cohen-Macaulay local ring,
Roberts obtained an $R$-module $M$ of finite length and finite projective
dimension with the property above.

In 2000, A.\ Singh and the author \cite{MS} used similar techniques
to obtain such a module over a five-dimensional Gorenstein local ring.
This shows that Theorem \ref{Duttacmgor} does not extend to all Gorenstein
local rings.

In fact, in each example above, the module $M$ satisfies an inequality
\[
\ell(F^n(M)) < \ell(M) p^{nd}
\] 
for all $n\gg 0$, in contrast to
modules of finite length over a complete intersection ring, which {\it
all}\/ satisfy the inequality $\ell(F^n(M)) \geq \ell(M) p^{nd}$, cf.\
Theorem \ref{Duttaci} (i).  
\end{exa}

\section{Complete intersection rings}\label{cisection}

In this section we discuss a class of rings for which the best
possible results characterizing finite projective dimension in terms
of homological (cf.\ \S\ref{rigidity}) and numerical (cf.\
\S\ref{cinumerical}) conditions involving the Frobenius endomorphism
hold. This case was rendered more approachable by techniques
particular to this class of rings. These are given in Section
\ref{tools}, before the proofs, which are gathered in Section
\ref{proofs}.

To introduce the next notion we recall that by Cohen's Structure
Theorem, cf.\ \S\ref{kunzsection}, every complete local ring has a
{\em Cohen presentation\/}, that is, a presentation as a homomorphic
image of a regular local ring.

\begin{ssdefn*} 
A local ring $R$ is said to be {\em complete intersection\/} if in
some (equivalently, in every) Cohen presentation of its $\fm$-adic
completion $\wh R$, the defining ideal is generated by a regular
sequence.  A Noetherian ring $R$ is said to be {\em complete
intersection\/} if for every maximal ideal $\fm$ the local ring
$R_\fm$ is complete intersection.  
\end{ssdefn*}

\subsection{Homological conditions for finite projective dimension}
\label{rigidity}

The results in Section \ref{kohlee}, especially Theorem \ref{KL},
raise the question whether the ultimate generalization of Theorem
\ref{H} holds: If {$\Tor_i^R(M,\up{\varphi^n}R)=0$} for some $i>0$ and
some $n>0$, then $M$ has finite projective dimension. A positive
answer is known over complete intersection rings.

\begin{ssthm}[Avramov-Miller, {\cite[Main Thm.]{AM}}]
\label{AvMi}
Let $R$ be a complete intersection ring, and 
$M$ a (possibly infinitely generated) $R$-module.

If $\Tor_i^R(M,\up{\varphi^n}R)=0$ for some fixed $i,n>0$ then
$\Tor_j^R(M,\up{\varphi^n}R)=0$ for all $j\ge i$; if, furthermore, $M$ is
finitely generated, then $\projdim M<\infty$, and hence 
$\Tor_j^R(M,\up{\varphi^n}R)=0$ for all $j>0$. 
\end{ssthm}

The first statement says that $\up{\varphi^n}R$ is rigid 
for any complete intersection ring $R$. 

It suffices to treat the case when the ring $R$ is local.  By
\ref{completion} below one may assume that $R$ is complete. The proof
in the complete local case exploits a very concrete factorization of
the $n$-th Frobenius endomorphism $\varphi^n$ as a flat map followed
by a surjection, discussed in Section \ref{tools}, and the notion of
complexity.  The factorization allows the possibly infinitely
generated module $\up{\varphi^n}R$ to be replaced by a cyclic module
over a (different) complete intersection ring. Furthermore, the
asymptotic homological behavior of this cyclic module is transparent
due to an explicit resolution constructed by Tate \cite{T}.  Dutta
\cite{D2} has given a variation on the proof that uses the
factorization, but avoids the notion of complexity.

In Section \ref{proofs} we sketch the proof from
\cite{AM}, incorporating a simplification of the first part
given in \cite{D2}.

\subsection{Numerical conditions for finite projective dimension}
\label{cinumerical}

The numerical functions discussed in Section \ref{numericalgens} 
behave well for complete intersection rings. 

\begin{ssthm}[Dutta, {\cite[Thm.~1.9]{D1}}]
\label{Duttaci}
Let $R$ be a complete intersection local ring, and 
let $M$ an $R$-module of finite length. Then 
\begin{enumerate}
\item[\rm{(i)}]
$\ell(F^n(M)) \geq \ell(M) p^{nd}$ for all $n\geq 0$, and, 
\item[\rm{(ii)}]  $\ell(F^n(M)) = \ell(M) p^{nd}$ for all $n\geq 0$ 
if $M$ has finite projective dimension.
\end{enumerate}
\end{ssthm}

Combining the original proof in \cite{D1} with ideas from \cite{AM},
and \cite{D2}, we prove this result in Section
\ref{proofs}. Local Chern character theory, described in Section
\ref{chern}, provides another proof.

Using the ideas from all three papers, \cite{D1}, \cite{AM}, and
\cite{D2}, we can derive a converse of Theorem \ref{Duttaci}.
The proof can be found in Section \ref{proofs}.

\begin{ssthm}
\label{length} 
Let $R$ be a complete intersection local ring, and let
$M$ be an $R$--module of finite length. 
If $\ell(F^n(M)) = \ell(M) p^{nd}$ for some $n>0$, 
then $M$ has finite projective dimension. 
\end{ssthm}

The condition on lengths is rendered more concrete 
if we recall that for a module $M$ with presentation 
\[
R^s \xra{\ [a_{ij}]\ } R^t \lra M \lra 0,
\]
the module $F^n(M)$ is simply the cokernel of 
the matrix $[a_{ij}^{\ p^n}]$.
In particular, we obtain a quick way to check the finiteness of
projective dimension of $\fm$-primary ideals.

\begin{sscor}
Let $I$ be an $\fm$-primary ideal in a complete intersection local ring $R$. 
Then $I$ has finite projective dimension if and only if 
$\ell(R/I^{[p]}) = p^d \ell(R/I)$.
\qed
\end{sscor}

\subsection{Main Ingredients}
\label{tools}

This section includes the basic tools that we need for the proofs of
the results above. It is these constructions that make the Frobenius
endomorphism in the case of complete intersection rings so amenable to
study.

We begin with a factorization, defined for any local ring $R$
of characteristic $p$; it is particularly useful in the case that $R$
is complete intersection. We will use

\begin{ssntn*}
For any homomorphism $\alpha\colon A\to B$ of
commutative rings, we let $\up{\alpha}B$ denote the $A$-$B$-bimodule
$B$ with $A$ acting through $\alpha$ and $B$ acting through $\id_B$,
that is, $a\cdot b'=\alpha(a)b'$ and $b'\cdot b=b'b$ for all $a\in A$,
$b'\in\up\alpha B$, $b\in B$.
\end{ssntn*}

\begin{factor}
\label{factorization} 
Let $R$ be a local ring of characteristic $p$.  By Cohen's Structure
Theorem, cf.\ \S\ref{kunzsection}, the ring $\wh R$ is a
residue ring of a ring of formal power series $Q=k[[\bst]]$ on
indeterminates $\bst=t_1, \dots,t_e$, 
say $\wh R=Q/I$, with $I\subseteq \fm_Q^2$. Set
\[
S=\wh R\otimes_Q \up{\varphi^n}Q=Q/I^{[p^n]}.
\]

Let $\sigma\colon S\to\wh R$ be the canonical
surjection, and let $\rho$ denote the composition 
\[ 
R\xra{\ \iota\ }
\wh R=\wh R\otimes_Q Q\xra{\wh R\otimes_Q\varphi_Q^n} 
\wh R\otimes_Q\up{\varphi^n}Q=S\,.  
\] 
Then 
$\sigma\rho=\varphi_{\wh R}^n\iota=\iota\varphi_R^n$;
as $\iota$ and $\varphi^n\colon Q\to Q$
are local flat homomorphisms, so is $\rho$.

Since $\wh R$ is flat over $R$ and $\rho$ is flat, we have isomorphisms
\begin{equation}\label{changerings}
\Tor_i^R(M,\up{\varphi^n}R)\otimes_R \wh R\cong
\Tor_i^R(M,\up{\varphi^n\iota}\wh R)\cong
\Tor_i^S(M\otimes_R\up\rho S,\up\sigma\wh R).
\end{equation}

Concretely, the endomorphism $\varphi^n_{\wh R}$ factors as
\[
\begin{CD}
\wh R @>{\psi=\wh R\otimes_Q\varphi_Q^n}>> S @>{\sigma}>> \wh R \\
@| @| @| \\
Q/I @>{}>>
Q/I^{[p^n]} @>{}>>Q/I
\end{CD}
\]
where $\psi(\overline{r})=\overline{r^{p^n}}$ and $\sigma(s)=\overline{s}$.
If $I=(x_1, \dots, x_r)$, 
then $I^{[p^n]}=(x_1^{p^n}, \dots, x_r^{p^n})$.  
In particular, if $R$ is complete intersection, then the ring $S$ is 
also complete intersection and has the same codimension as $R$.  
\end{factor}

The following remark is useful for reducing to the complete case. 

\begin{compl} 
\label{completion} 
If $\iota\colon R\to\wh R$ denotes the
canonical map into the $\fm$-adic completion, then 
\[
\Tor_i^R(M,\up{\varphi^n\iota}\wh R)\cong
\Tor_i^R(M,\up{\varphi^n}R)\otimes_R \wh R\cong 
\Tor_i^{\wh R}(M\otimes_R \wh R,\up{\varphi^n}\wh R)
\]
for all $i\geq 0$ since $\iota\varphi^n_R=\varphi^n_{\wh R}\iota$. 
For $i=0$ this yields
\[
F^n_R(M)\otimes_R \wh R\cong 
F^n_{\wh R}(M\otimes_R \wh R)
\]
and thus an equality of lengths 
$\ell_R(F^n_R(M))=\ell_{\wh R}(F^n_{\wh R}(M\otimes_R \wh R))$. 
Since $\wh R$ is faithfully flat over $R$, 
these facts enable us to assume that $R$ complete in the proofs below. 
\end{compl} 

We describe next a filtration originating in \cite{D1}
that conveniently complements the factorization above. 

\begin{filtr} 
\label{filtration} 
Assume now that $R$ is complete. 
By Cohen's Structure Theorem, cf.\ \S\ref{kunzsection}, 
$R$ can then be written as a residue ring of a regular ring $Q$ by
a $Q$-regular sequence $\bsx=x_1,\dots, x_c\in\fm^2$, that is, $R=Q/(\bsx)$.
We use the notation from the factorization in 
Remark \ref{factorization}. 

The $Q$-module $S=Q/(\bsx^{p^n})$ has a filtration 
\[
0=S_{p^{nc}}\subset S_{p^{nc}-1}\subset\cdots\subset S_{1}\subset S_{0}=S
\]
with subquotients isomorphic to $R$; it produces exact sequences
\[
0 \lra S_{k+1} \xra{\ \tau_k\ } S_{k} \xra{\sigma_{k}}R \lra 0
\quad\text{for}\quad k=0, \dots, p^{nc}-1\,
\]
with $\sigma_{0}$ equal to the map $\sigma$ in the factorization in 
Remark \ref{factorization}. 
Tensoring each sequence on the left with the $S$ module 
$M'=M\otimes_R\up\rho S$,  
we obtain, for $k=0, \dots, p^{nc}-1$, long exact sequences
\begin{equation}\label{les}
\lra \Tor_1^S(M',\up\sigma R)  
\xra{\delta_k} M'\otimes_S S_{k+1} \xra{\ \tau_k\ } 
M'\otimes_S S_{k} \xra{\sigma_{k}}
M'\otimes_S \up\sigma R \lra 0
\end{equation}
In each long exact sequence, the last term can be recognized as 
\begin{equation}\label{fnrm}
M'\otimes_S \up\sigma R=
(M\otimes_R\up\rho S)\otimes_S R\cong M\otimes_R\up{\varphi^n} R 
=F^n_R(M).
\end{equation}
Similarly, in the sequence for $k=0$, since $S_0=S$ 
the second-to-last term is 
\begin{equation}\label{fnqm}
M'\otimes_S S_0\cong M'=
M\otimes_R\up\rho S \cong M\otimes_Q\up{\varphi^n} Q = F^n_Q(M).
\end{equation}
\end{filtr} 

\subsection{Proofs}
\label{proofs}

We are now ready to give the proofs of the results in \ref{rigidity}
and \ref{cinumerical}. In each proof, we use
the notation and constructions from Remarks \ref{factorization} and
\ref{filtration}. In particular, once we have reduced to the situation where
$R$ is complete, setting $M'=M\otimes_R\up\rho S$, 
we use repeatedly the isomorphisms (\ref{changerings}) from \ref{factorization}
\[
\Tor_j^R(M,\up{\varphi^n}R)\cong\Tor_j^S(M',\up\sigma R).
\]

The first proof is from \cite{AM} but is modified to include the 
simplification of the first part given in \cite{D2}.

\begin{proof}[Sketch of Proof of Theorem \ref{AvMi}]
By \ref{completion} we may assume that $R$ is complete. 

The first part is proved by induction on $j\geq i$.  Suppose that
$\Tor_j^R(M,\up{\varphi^n}R)=0$ and thus {$\Tor_j^S(M',\up\sigma R)=0$}.
Since $S_{p^{nc}}=0$, the long exact sequences
(\ref{les}), beginning with the one for 
$k=p^{nc}-1$ and ending with the one for $k=1$, yield
{$\Tor_j^S(M',S_1)=0$}. Since $S_0=S$, it is automatic that
{$\Tor_{j+1}^S(M',S_0)=0$} and so the long exact sequence for $k=0$ 
yields {$\Tor_{j+1}^S(M',\up\sigma R)=0$}, as desired.

For the second part, by passing to a syzygy module of $M$ we may
assume $i=1$, and so {$\Tor_j^R(M,\up{\varphi^n}R)=0$} for all $j\geq
0$. Equivalently, {$\Tor_j^S(M',\up\sigma R)=0$} for all $j\geq 0$,
and so it follows from \cite[Prop.~2.1]{Mi1} that
\[
\cxy_S M' +\cxy_S \up\sigma R = \cxy_S (M'\otimes_S\up\sigma R)
\]
where $\cxy_S$ denotes complexity, cf.\ Appendix \ref{alphacx}. 
The right-hand side is at most equal to $c=\codim S$, cf.\ Appendix
\ref{alphacx}.  On the other hand, by a result of Tate
\cite[Thm.~6]{T}, an explicit resolution of $\up\sigma R$ over the
complete intersection ring $S$ is known, and it yields $\cxy_S
\up\sigma R=c$.  This forces $\cxy_S M'=0$ or, equivalently,
$\pd_SM'<\infty$. Since $\rho$ is faithfully flat, this implies that 
$\pd_RM<\infty$. 
\end{proof}

Combining the original proof in \cite{D1} with ideas from \cite{AM},
and \cite{D2} yields

\begin{proof}[Proof of \ref{Duttaci}]
By \ref{completion} we may assume that $R$ is complete. 

Keeping the isomorphisms (\ref{fnrm}) and 
(\ref{fnqm}) in mind, we see by counting lengths in the exact 
sequences (\ref{les}) that
\[
\ell(F^n_Q(M))\leq \ell(F^n_R(M))p^{nc},
\]
with equality holding if and only if $\delta_j=0$ for $j=0, \dots, p^{nc}-1$. 
Since $Q$ is regular of dimension $d+c$, the left-hand side equals 
$\ell(M)p^{n(d+c)}$ by Remark \ref{lengthsregular}; thus, 
\[
\ell(F^n_R(M))\geq \ell(M)p^{nd},
\]
with equality holding if and only if $\delta_j=0$ for $j=0, \dots, p^{nc}-1$. 

If $M$ has finite projective dimension, then
$\Tor_i^R(M,\up{\varphi^n}R)=0$ for all $i>0$ by Theorem \ref{PS} and
thus $\Tor_i^S(M',\up\sigma R)=0$ for all $i>0$. Therefore from the 
exact sequences (\ref{les}) we get $\delta_j=0$ for all $j$, giving the
conclusion of (ii).
\end{proof}

The remaining result is now almost immediate. 

\begin{proof}[Proof of \ref{length}]
As shown in the preceding proof, the equality 
$\ell(F^n_R(M)) = \ell(M) p^{nd}$ implies that $\delta_j=0$ for all $j$.
In particular, $\delta_0$ is the zero map, and the sequence (\ref{les}) 
yields {$\Tor_1^S(M',\up\sigma R)=0$} since {$\Tor_1^S(M',S)=0.$}
Therefore, one obtains $\Tor_1^R(M,\up{\varphi^n}R)=0$, 
so $M$ has finite projective dimension by Theorem \ref{AvMi}.
\end{proof}

\subsection{Growth of Tors}
\label{growthtors}

Theorem \ref{extremal} yields a comparison of the asymptotic
homological behaviors of the modules $\up{\varphi^n}R$ and $k$; more
specifically it relates the growths of the sequence
$\ell_{\varphi^n}\big(\Tor_i^R(k,\up{\varphi^n}R)\big)$ to the Betti
numbers $\ell_R\big(\Tor_i^R(k,k)\big)$ of the residue field.  In this
section we discuss an extension of this for complete intersection
rings in which $k$ is replaced by any $R$-module $M$ of finite length.

\begin{ssrem}\label{limitAM2}
The Betti numbers of $M$ have quasi-polynomial behavior: 
by \cite[Cor.~4.1]{G} and \cite[Thm.~4.1]{Av1} there exist
polynomials $b_\pm(t)\in\BQ[t]$ such that
\begin{gather*}
\ell_R\big(\Tor^R_i(M,k)\big)=
\begin{cases}
b_+(i) &\text{for all}\quad i=2s\gg0\,;\\
b_-(i) &\text{for all}\quad i=2s+1\gg0\,;
\end{cases}\\
\deg b_+(t)=\deg b_-(t)=\cxy_R M-1\,.
\end{gather*}
\end{ssrem}

Similarly, we ask about the rate of growth of
$\ell_{\varphi^n}(\Tor_i^R(M,\up{\varphi^n}R))$, where the length is measured
as in Appendix \ref{fgsection}.  The next result
shows that if $\ell_R(M)$ is finite and $n$ is fixed, the numbers
$\ell_{\varphi^n}(\Tor_i^R(M,\up{\varphi^n}R))$ grow at the same rate as do the
numbers $\ell_R(\Tor_i^R(M,k))$.

\begin{ssthm}[Avramov-Miller, {\cite[Main Thm.]{AM}}] 
\label{AM2}
Let $R$ be a complete intersection local ring and $M$ an $R$-module. 
If $\ell_R(M)<\infty$, then for every $n>0$ 
there exist polynomials $h_\pm(t)\in\BQ[t]$ such that
\begin{gather*}
\ell_{\varphi^n}\big(\Tor_i^R(M,\up{\varphi^n}R)\big)=
\begin{cases}
h_+(i) &\text{for all}\quad i=2s\gg0\,;\\
h_-(i) &\text{for all}\quad i=2s+1\gg0\,;
\end{cases}\\
\max\{\deg h_+(t)\,,\,\deg h_-(t)\}=\cxy_R M-1\,.
\end{gather*}
\end{ssthm}

The proof involves using Matlis 
duality to convert to Ext modules, where an algebra structure and the theory 
of the complexity of pairs of modules is available. 

\section{Asymptotic numerical conditions}\label{asymptotic}

\subsection{Asymptotic criteria}\label{asymptoticcriteria}

Kunz \cite{K}, cf.\ \S\ref{kunzsection}, introduced the numerical function
$n\mapsto \ell(R/\fm^{[p^n]})$ as an analog of the Hilbert-Samuel function
$n\mapsto \ell(R/\fm^n)$.

More generally, Monsky in 1983 considered the function $n\mapsto
\ell(M/I^{[p^n]}M)$ for any $\fm$-primary ideal $I$ and finitely generated
$R$-module $M$ and named it the {\em Hilbert-Kunz function\/} of $M$.
As an analog of the fact that the Hilbert-Samuel function is a polynomial
of degree $d$ for $n\gg 0$ with leading coefficient $e(\fm,R)$ defining
the multiplicity of $R$, he proved:

\begin{ssthm}[Monsky, {\cite[Thm.~1.8]{Mo}}]
\label{monsky}
Let $R$ be a local ring of characteristic $p$, $I$ an $\fm$-primary 
ideal of $R$, and $M$ a finitely generated $R$-module. There exists 
a constant $C$ such that 
\[
\ell(M/I^{[p^n]}M)=C p^{nd} + O(p^{n(d-1)}),
\]
where $O(p^{n(d-1)})$ represents a term whose absolute value is less 
than $Bp^{n(d-1)}$ for some constant $B$ and all $n\gg 0$.
\end{ssthm}

In the case that $M=R$, we can rewrite this function in terms of the
Frobenius functors as $\ell(F^n(R/I))$ and it is usually called the {\em
Hilbert-Kunz function\/} of $I$. By Theorem \ref{monsky}, 
the limit
\[
\lim_{n\to\infty} \frac{\ell(F^n(R/I))}{p^{nd}}, 
\]
then exists; it is known as the {\em Hilbert-Kunz multiplicity\/} 
$e_{\rm HK}(I,R)$ of $I$. If $I=\fm$, the limit is called the 
Hilbert-Kunz multiplicity of the local ring $R$. 

The Hilbert-Kunz function and multiplicity have been studied by a
number of authors, beginning with Monsky in \cite{Mo} and Han and
Monsky in \cite{HM}, and tend to be very complicated when the ideal in
question does not have finite projective dimension. Defined similarly
to the usual multiplicity of an ideal $I$, the Hilbert-Kunz
multiplicity however seems to reflect arithmetic
information about the ideal. Not much is known yet about its
interpretations or even whether it is always a rational number.

Recently Watanabe and Yoshida characterized regularity of a local ring
in terms of an asymptotic version of condition (c) in Kunz's Theorem.
A simplified proof was given later by Huneke-Yao \cite{HY}. This
result is analogous to the classical characterization of regularity in
terms of the Hilbert-Samuel multiplicity that says that if $R$ is
unmixed, then it is regular if and only if $e(\fm,R)=1$. 
Recall that $R$ is unmixed if $\dim \wh R= \dim \wh R/\fp$ 
for every associated prime ideal $\fp$ of $\wh R$. 

\begin{ssthm}[Watanabe-Yoshida, {\cite[Thm.~1.5]{WY}}]
\label{WY}
Let $R$ be a unmixed local ring 
of characteristic $p$. 
Then the following conditions are equivalent.
\begin{enumerate}
\item[\rm{(a)}]
$R$ is regular.
\item[\rm{(c$^*$)}]
$e_{\rm HK}(\fm,R)=1$   
\end{enumerate}
\end{ssthm} 

\begin{ssrem}
The condition that $R$ is unmixed is necessary as 
the following well-known example shows: 
Let $k$ be a field of characteristic $p$, and 
$R$ the ring $k[[x,y]]/(x^2,xy)$ with maximal ideal $\fm=(x,y)$. 
Then $e_{\rm HK}(\fm,R)=1$ and yet $R$ is not regular. 
\end{ssrem}

Dutta extended Hilbert-Kunz functions of ideals to a module setting  
by using Frobenius functors. As in the previous section, let $M$
denote a module of finite length. Seibert \cite{Se} proved that 
there is a constant $C$ such that 
\[
\ell(F^n(M))=C p^{nd} + O(p^{n(d-1)}),
\]
so that the limit
\[
\lim_{n\to\infty} \frac{\ell(F^n(M))}{p^{nd}} 
\]
exists. Theorem \ref{Duttaci} implies that if $M$ is a module of
finite projective dimension over a complete intersection ring, then
this limit equals $\ell(M)$. Furthermore, over complete intersection
rings the converse was proved in \cite{Mi}, yielding the corresponding
module version of Theorem \ref{WY} over these rings. We give the proof
here since it can be shortened substantially now in view of the more
recent results.

\begin{ssthm}[Miller, {\cite[Thm.~2.1]{Mi}}]
\label{limitnumerical}
Let $R$ be a complete intersection local ring of characteristic $p$, and let
$M$ be an $R$--module of finite length. 
The following conditions are then equivalent. 
\begin{enumerate}
\item[\rm{(a)}]
$M$ has finite projective dimension;
\item[\rm{(c$^*$)}]
$\displaystyle{\lim_{n\to\infty}
\frac{\ell(F^n(M))}{p^{nd}}=\ell(M).}$
\end{enumerate}
\end{ssthm}

\begin{proof}
By Theorem \ref{Duttaci} with $n=1$, applied to the modules 
$F^i(M)$ for $i\geq 0$, 
we see that the sequence 
\[
\left ( 
\frac{\ell(F^n(M))}{p^{nd}}
\right )_{n=0}^{\infty}
\]
is nondecreasing. 
If $\projdim(M)<\infty$, then the limit equals $\ell(M)$ by 
Theorem \ref{Duttaci}.  

Conversely, if the limit of the sequence equals $\ell(M)$, namely, its
initial term, then all terms are equal, i.e., 
$\ell(F^n(M)) = \ell(M) p^{nd}$ for all $n\geq 0$,
and so $M$ has finite projective dimension by Theorem \ref{length}. 
\end{proof}


\begin{sscor}
Let $R$ be a complete intersection local ring of characteristic $p$. 
If $I$ is an $\fm$-primary ideal, then it has finite projective 
dimension if and only if 
\begin{xxalignat}{3}
&{\hphantom{\square}}& e_{\rm HK}(I,R)&=\ell(R/I)  
&&\square
\end{xxalignat}
\end{sscor}

\subsection{Local Chern characters}\label{chern}

We cannot expect results like Theorem \ref{limitnumerical}
to hold over arbitrary local rings since, in fact, 
the modules in Example \ref{exs} satisfy 
\[
\lim_{n\to \infty}\frac{\ell(F^n(M))}{p^{nd}} \neq \ell(M).
\]
Why such examples can exist over non-complete intersection
rings is explained by the theory of local Chern characters of Baum,
Fulton, and MacPherson, \cite{BFM} and \cite{Fu}. It was developed
further by P.\ Roberts in \cite{Ro1} for use in his proof of Serre's
Vanishing Conjecture for rings of mixed characteristic, where his
result on the commutativity of the local Chern characters was the main
ingredient.  He used the theory again in an essential way in his proof
of the New Intersection Theorem in mixed characteristic
\cite{RobMSRI}.

Following the exposition in \cite{RobMSRI} and in \cite{RobBook},
we briefly sketch 
the explanation given by Roberts of the obstruction to an equality
\[
\lim_{n\to \infty}\frac{\ell(F^n(M))}{p^{nd}} = \ell(M)
\]
for modules $M$ of finite projective dimension.

Assume first that $R$ is an arbitrary complete local ring.  
For a closed subset $X$ of $\Spec R$ let
$A_*(X)_{\mathbb Q}=\bigoplus A_i(X)_{\mathbb Q}$ denote the rational
Chow group of $X$, that is, $A_i(X)_{\mathbb Q}$ is the $\mathbb
Q$-vector space on cycles of codimension $i$ modulo rational equivalence.
Given a bounded complex $L_\cx$ of finitely generated free $R$ modules
with finite length homology, there exists a local Chern character
\[
\ch(L_\cx) = \ch_d(L_\cx) + \ch_{d-1}(L_\cx) + \cdots + \ch_0(L_\cx)
\]
where, for each $i$, 
\[
\ch_i(L_\cx)\colon A_i(\Spec R)_{\mathbb Q} \to 
A_{0}(\{ \fm \} )_{\mathbb Q}\cong \mathbb Q.
\]
If $L_\cx$ is the resolution of a module $M$ of finite length and finite
projective dimension, then by a special case of the Local Riemann-Roch 
Formula (cf.\ \cite[Thm.~12.6.1]{RobBook}) there exists an element
$$
\tau(R) = \tau_d(R) + \cdots + \tau_0(R) \in A_*(\Spec R)_{\mathbb Q},
$$ 
called the Todd class of $R$, such that 
$$
\ell(M) = \ch(L_\cx)(\tau(R)) = \sum_{i=0}^d \ch_i(L_\cx)(\tau_i(R)).
$$

Now suppose that $R$ has positive characteristic $p$ and perfect
residue field and that the Frobenius endomorphism is a finite map.
The following result follows from the facts that 
$F^n(L_\cx)$ is a resolution of $F^n(M)$ (cf.\ Theorem \ref{PS}) and 
that local Chern characters are compatible with finite maps. 

\begin{ssthm}[P.\ Roberts, {\cite[Prop.\ 12.7.1]{RobBook}}]
\label{chernd} 
With the hypotheses above, 
\[
\ell(F^n(M))=\ch(F^n(L_\cx))(\tau(R)) = 
\sum_{i=0}^d p^{ni}\ch_i(L_\cx)(\tau_i(R)),
\]
and thus 
\[
\lim_{n\to\infty}\frac{\ell(F^n(M))}{p^{nd}} = \ch_d(L_\cx)(\tau_d(R)). 
\]
\end{ssthm}

This result, as well as Roberts' proof of the Vanishing Theorem, 
are motivations for the following definition.

\begin{ssdefn}[Kurano, {\cite[Def.\ 2.1]{Kurobring}}]
Let $R$ be an algebra of finite type over some regular local ring. 
The ring $R$ is {\em Roberts\/} if $\tau_i(R)=0$ for all $i<d$.
\end{ssdefn}

\begin{ssdefn}[Kurano, {\cite{Kunumer}}]
Let $R$ be an algebra of finite type over some regular local ring. 
The ring $R$ is {\em numerically Roberts\/} 
if for any finite free complex $L_\cx$ with finite length homology 
$\ch_i(L_\cx)(\tau_i(R))=0$ for all $i<d$. 
\end{ssdefn}

\begin{sscor}\label{=numrob}
If $R$ is a numerically Roberts ring and $M$ a module of finite projective 
dimension, then $\ell(F^n(M)) = \ell(M) p^{nd}$ for all $n\geq 0$.
\qed
\end{sscor}
    
Certain numerically Roberts rings have identified by using
the following result.

\begin{ssthm}[P.\ Roberts, 
{\cite[Thm.~2]{RobMacRae}, \cite[Thm.~12.4.4]{RobBook}}]
\label{tau}
Let $R$ be a complete local ring of characteristic $p$. 
\begin{enumerate}
\item[\rm{(i)}]
If $R$ is complete intersection, then $\tau_i(R)=0$ if $i<d$, i.e., 
$R$ is Roberts. 
\item[\rm{(ii)}]
If $R$ is Gorenstein, then $\tau_{d-i}(R)=0$ if $i$ is odd. 
\item[\rm{(iii)}]
If $\dim R > 0$, then $A_0(\Spec R)_{\mathbb Q}=0$. 
\item[\rm{(iv)}]
If $L_\cx$ has finite length homology and $\dim R > 1$, 
then $\ch_1(L_\cx)=0$. 
\end{enumerate}
\end{ssthm}
 
For lack of a suitable reference, we briefly explain (iii): letting 
$\fp$ be a submaximal prime ideal of $R$ and $x$ a non-unit outside 
of $\fp$, so that $\fm$ is the only prime ideal in the support of  
$R/((x)+\fp)$, we see that $\divcl(x,\fp)=\ell(R/((x)+\fp))[R/\fm]$ 
in the Chow group of $R$. Thus $[R/\fm]$ is zero in the rational Chow group. 

Theorems \ref{chernd} and \ref{tau} reaffirm that the equality
$\ell(F^n(M)) = \ell(M) p^{nd}$ holds for modules $M$ of finite
projective dimension over the rings in each of the classes in Theorems
\ref{Duttacmgor} and \ref{Duttaci} since in each case
$R$ is numerically Roberts.
However, in general the terms $\ch_i(L_\cx)(\tau_i(R))$ for $i<d$ 
account for a possible inequality 
\[
\lim_{n\to \infty}\frac{\ell(F^n(M))}{p^{nd}} \neq \ell(M),
\]
such as in Roberts' example described in Example \ref{exs}. In 1996
Kurano \cite{Ku} found a five-dimensional Gorenstein ring that is not
Roberts ($\tau_3(R)\neq 0$). The example of Singh and the author
\cite{MS}, cf.\ Example \ref{exs}, provided the first 
five-dimensional Gorenstein ring that is not numerically Roberts.

Very recently, Roberts and Srinivas \cite{RS}\! used deep K-theoretic
results of Thomason and Trobaugh to construct large families of rings
which are not numerically Roberts: these include affine cones of
$\mathbb P^n \times \mathbb P^m$ studied by Kurano in \cite{Ku}. More
generally, Roberts and Srinivas consider affine cones of smooth
varieties which have a nondegenerate intersection pairing, and obtain
conditions under which Roberts and numerically Roberts are equivalent
notions. More precisely, they show that for such rings whenever
$\tau_i(R)\neq 0$ for some $i<d$, there is a complex $L_\cx$ with
$\ch_i(L_\cx)(\tau_i(R))\neq 0$. Whether $\tau_i(R)$ vanishes can be
determined, modulo Grothen\-dieck's {\em Standard Conjectures\/}, by examining
the related theory of topological Todd classes. The case of
Grassmanians was later studied by Kurano and Singh \cite{KS}, who
determine which of these have Roberts homogeneous coordinate rings.

\section{Questions}\label{questions}

In this section we consider some natural questions arising from the
results discussed in this survey.

The question remains whether $\up{\varphi^n}R$ is a ``test module"
for finite projective dimension over any local ring:

\begin{conj}\label{conj}
 Let $M$ be an $R$-module. 
If $\Tor_i^R(M,\up{\varphi^n}R) = 0$ for some fixed $i,n>0$, 
then $M$ has finite projective dimension. 
\end{conj}

The evidence seems strong: It holds for complete intersection rings by
Theorem \ref{AvMi}. At the other extreme, it holds also for rings of
depth zero as long as $n\geq\log_p(\crs_p(R))$ by Corollary
\ref{KLdepth0} and for Cohen-Macaulay rings of dimension one as long
as $n\geq\log_p(\drs_p(R))$ by Corollary \ref{KLdepth1}. 

As discussed in Section \ref{numericalgens}, if $M$ is a finitely
generated $R$-module of finite projective dimension, then an equality
$\ell(F^n(M)) = \ell(M) p^{nd}$ holds in the graded case, cf.\ Theorem
\ref{PSgraded}; it also holds for complete intersection rings, cf.\
Theorem \ref{Duttaci} and for some classes of rings of small
dimension, cf.\ \ref{Duttacmgor} and \ref{Duttaci} (indeed, it holds over
any ring that is numerically Roberts, cf.\ \ref{=numrob}).  However, by
Example \ref{exs}, it may fail in general. When it does hold,
there is a possibility for a converse; we list the graded case,
perhaps the most approachable.

\begin{ques}
Let $R=\bigoplus_{i\geq 0}R_i$ be a graded Noetherian ring of 
dimension $d$ with $R_0$
Artinian, and let $M$ be a graded $R$-module of finite length.  
If the equality $\ell(F^n(M)) = \ell(M) p^{nd}$ holds for all $n\geq 0$, does
$M$ then have finite projective dimension? 
\end{ques}

Next we consider, for a finitely generated module $M$ over a local
ring $R$, the properties of $F^n_R(M)$. By Remark \ref{dp-exs},
whenever $M$ has finite projective dimension, the modules $F^n_R(M)$
and $M$ have the same depth, but this fails in general.  On the other
hand, it is possible for the sequence $\{ \depth F^n_R(M)\}_n$ to be 
constant even when $\projdim_RM$ is infinite: modules of finite 
length give such examples.

\begin{ques} 
Is there a stable value of $\depth F^n_R(M)$ as $n\to\infty$?  
For which $R$-modules $M$ is the depth of $F^n_R(M)$ constant for all $n$?  
\end{ques}

Finally, we consider asymptotic homological properties of the Frobenius
endomorphism. Let $\mu_{\varphi^n}(-)$ denote the minimal number of
generators.  A natural question arising from Theorem \ref{AM2} is:

\begin{ques} 
If $R$ is a complete intersection local ring and $M$ a finitely
generated $R$-module, does then the function $i\mapsto 
\mu_{\varphi^n}\big(\Tor_i^R(M,\up{\varphi^n}R)\big)$ have quasi-polynomial
behavior for $i\gg 0$? 
\end{ques}

If $\ell_R(M)<\infty$, this question has an affirmative answer by the
proof of Theorem \ref{AM2}, cf.\ \cite[\S2]{AM} for details.

It is known that $\ell_R\big(\Tor_i^R(M,k)\big)$ has maximal growth
when $M$ is equal to the residue field, in the following precise sense: for
any finitely generated $R$-module $M$, there exists a constant $C_M$
such that 
\[
\ell_R\big(\Tor_i^R(M,k)\big)\leq C_M\cdot\ell_R\big(\Tor_i^R(k,k)\big)
\, {\textrm{ for all }}i\gg0.
\]
In terms of complexity and curvature,
defined in Appendix \ref{alphacx}, this implies 
\[
\cxy_R(M)\leq\cxy_R(k) \qquad {\textrm{and}} \qquad  \curv_R(M)\leq\curv_R(k).
\]
This suggests the following question for Tors against 
the powers of the Frobenius endomorphism. 

\begin{ques} 
For any finitely generated module $M$ over a local ring $R$, 
does there exist a constant $C_M$ such that an inequality 
\[
\mu_{\varphi^n}\big(\Tor_i^R(M,\up{\varphi^n}R)\big)\leq 
C_M\cdot\ell_{\varphi^n}\big(\Tor_i^R(k,\up{\varphi^n}R)\big),
\]
holds for all $i\gg0$?
\end{ques}

\begin{appendix}
\section{Module structures}\label{modulestructure}

The power of methods based on the Frobenius endomorphism is partly due
to the fact that it provides every $R$-module with two different actions of
$R$. This creates different $R$-module structures on tensor products,
Tors, etc., which may be confusing for first-time users. In this
appendix we review these structures and the basic properties they inherit
from the original $R$-module.

In order to differentiate between the various actions of $R$, we first
review the necessary facts in a more general\footnote{In
Appendix \ref{frobapp} we specialize to the case where $\alpha$ is the
$n$-fold composition  $\varphi^n\colon R\to R$ of the Frobenius
endomorphism of a Noetherian ring $R$, $M$ is an $R$-module, and
$N=R$.} set-up: $\alpha\colon A\to B$ is a homomorphism of commutative
Noetherian rings, $M$ is an $A$-module and $N$ is a $B$-module. We
always view $N$ as an $A$-module by restriction of scalars via
$\alpha$.

\subsection{Module structures on the Tors}\label{B-structure}

The tensor product $M\otimes_A N$ has a canonical structure of
a $B$-module\footnote{It is perhaps useful to think of $A$ and $B$ as
noncommutative, so the $A$-module structures on $M$ and $N$ have been
``used up" in forming the tensor product and only a $B$-structure
remains; it is via the $B$-structure that modules retain
properties such as finite generation and finite length.} via $N$,
given by $b(m\otimes n)=m\otimes nb$.

Similarly, each $\Tor_i^A(M,N)$ has a canonical structure of
$B$-module, via $N$. It can be defined by writing $\Tor_i^A(M,N)$ as
the homology of the complex $(L_\cx\otimes_A N)$, where $L_\cx$ is a
resolution of $M$ by free $A$-modules.

\subsection{Generation and length}\label{fgsection}

Some properties of $M$ are preserved by the Tors. 

\begin{fg} \label{fg}
If $M$ is finitely generated over $A$ and $N$ is finitely
generated over $B$, then $\Tor_i^A(M,N)$ is finitely generated over 
$B$. 
\end{fg}

\begin{proof} 
Let $L_\cx$ be a resolution of $M$ by finitely generated free
$A$-modules. Each module in the complex $L_\cx\otimes_A N$ is then a
finite direct sum of copies of $N$ and hence it is a finitely
generated $B$-module. Since $B$ is Noetherian the homology modules are
finitely generated as well.  
\end{proof}

\begin{fl} \label{fl}
Suppose that $N$ is a finitely generated $B$-module such that 
$\ell_B(N/\alpha(\fm) N)<\infty$ for every maximal ideal 
$\fm$ of $A$. If $M$ has finite length over $A$, then $\Tor_i^A(M,N)$
has finite length over $B$. 
\end{fl}

We use $\ell_\alpha(\Tor_i^A(M,N))$ to denote the length of the $B$-module 
$\Tor_i^A(M,N)$.

\begin{ssrem}\label{finfiber}
The condition on the lengths of $N/\alpha(\fm) N$ holds, in
particular, if $\alpha$ is a map with finite closed fibers, that is, if
$\ell_B(B/\alpha(\fm) B)$ is finite for every maximal ideal $\fm$ of
$A$.  
\end{ssrem}

\begin{proof}
Since $\Tor_i^A(M,N)$ is a finitely generated $B$-module by \ref{fg}, 
it is enough to show that its $B$-support consists only of maximal ideals.
Let $\fq$ be a non-maximal prime ideal in $B$, 
and let $\fp$ be its contraction to $A$. 
There are isomorphisms 
\[
\Tor^A_i(M,N)_{\fq} \cong
\Tor^A_i(M,N)_{\fp}\otimes_{A_{\fp}} B_{\fq} \cong
\Tor^A_i(M_{\fp},N_{\fq}). 
\]
If $\fp$ is maximal, then the set 
$\Supp(N/\alpha(\fp) N)=\Supp(B/\alpha(\fp) B)\bigcap\Supp(N)$ 
consists of maximal ideals by hypothesis. In this case, $N_{\fq}=0$ 
since $\alpha(\fp)\subseteq\fq$. If $\fp$ is not maximal, then 
$M_{\fp}=0$. 
\end{proof}

\subsection{Tor as an $A$-module}\label{A-structure}

Since $A$ is commutative, each $\Tor_i^A(M,N)$ has also a structure
of $A$-module, induced by the action of $A$ on $M$.  This endows
$\Tor_i^A(M,N)$ with a natural structure of $A$-$B$-bimodule for each
$i\geq 0$, cf.\ \cite[V.7]{Mac}.  Now the $B$-module structure gives rise
to another $A$-module structure by restriction of scalars via $\alpha$.

\begin{comp} \label{comp} 
The two $A$-module structures defined above on $\Tor_i^A(M,N)$ coincide. 
\end{comp}

\begin{proof}
For $i=0$ the conclusion follows from the equalities
\[
a(m\otimes n)=am\otimes n=ma\otimes n=
m\otimes an=m\otimes n\alpha(a)=(m\otimes n)\alpha(a)
\]
for any $m\in M$, $n\in N$ and $a\in A$. 
Let $L_\cx$ be a projective resolution of $M$ over $A$. 
By the same argument, the two $A$-module structures on the 
complex $L_\cx\otimes_A N$, and thus on its homology modules 
$\Tor_i^A(M,N)$, agree. 
\end{proof}

Applying \ref{fg} to the identity map of $A$ we get:

\begin{fg}\label{fgA}
If $M$ and $N$ are finitely generated over $A$, then the $A$-module 
$\Tor_i^A(M,N)$ is finitely generated. 
If, furthermore, $\ell_A(M)<\infty$, 
then $\ell_A(\Tor_i^A(M,N))<\infty$.
\qed
\end{fg}

Now suppose that $\alpha$ is local and that $A$ and $B$ have residue fields
$k$ and $l$, respectively. The map $\alpha$ induces an injection 
of residue fields $k\hra\fel$. 
Let $W$ be a $B$-module that has finite length when viewed as an
$A$-module. 
A filtration of $W$ over $B$ with $\ell_B(W)$ subquotients
isomorphic to $\fel$, which has length $\ell_k(\fel)$ over $A$, yields
\[
\ell_k(\fel)\cdot \ell_B(W)=\ell_A(W)< \infty
\]
In particular, $\ell_B(W)$ is finite, and so is $\ell_k(\fel)$ if
$W\neq 0$.  By \ref{fgA} the preceding discussion applies to
$W=\Tor_i^A(M,N)$, if $\ell_A(M)<\infty$ and $N$ is finitely generated
over $A$, cf.\ \ref{fl}.

\begin{flcomp} \label{flcomp}
Suppose that $\alpha\colon A\to B$ is a local homomorphism of local rings
with residue fields $k$ and $\fel$, respectively. If
$\ell_A(M)<\infty$ and $N$ is finitely generated over $A$, then
\begin{xxalignat}{3}
&{\hphantom{\square}}& \ell_A(\Tor_i^A(M,N))&=
\ell_k(\fel)\cdot \ell_\alpha(\Tor_i^A(M,N)) <\infty
&&\square
\end{xxalignat}
\end{flcomp}

\section{Complexity and curvature over a homomorphism}\label{alphacx}
\setcounter{subsection}{1}
\setcounter{ssthm}{0}

We extend the notions of complexity and curvature of a finite module
over a local ring $A$ to a relative situation. 

Let $\alpha\colon A\to B$ be a local homomorphism of local rings with
residue fields $k$ and $\fel$, respectively, such that $\ell_B(k\otimes_A
B)<\infty$. Let $N$ be a finitely generated $B$-module.  Then by \ref{fl}
we have $\ell_\alpha(\Tor_i^A(k,N))<\infty$ for all $i\geq 0$.  So we may
make the following definition.

\begin{ssdefn} \label{cxcurvextn}
Define the {\em complexity\/} and {\em curvature\/} of $N$ over $\alpha$ as 
\begin{align*}
\cxy_\alpha\!N&=
\inf\{t\in\mathbb N_0\mid \ell_\alpha(\Tor_i^A(k,N) \leq \beta\, i^{t-1}, 
\text{ some }\beta\in {\mathbb{R}}, \text{ all }i\gg0\} \\
\curv_\alpha\!N&=
\limsup_i \sqrt[i]{\ell_\alpha(\Tor_i^A(k,N))}
\end{align*}
\end{ssdefn}

When $\alpha$ is the identity map $\id_A$, 
we set 
\[
\cxy_AN=\cxy_{\,\id_A}\!N \qquad \textrm{and} \qquad
\curv_AN=\curv_{\,\id_A}\!N.
\] 
These ``absolute" notions have been used earlier, cf.\ \cite{Av2}.  We
review some standard facts for these: Clearly, an
$A$-module $N$ has finite projective dimension if and only if
$\cxy_AN=0$.  If $A$ is complete intersection, then
$\cxy_AN\leq\embdim A-\dim A$ by \cite[Cor.~4.1]{G}, with equality
when $N=k$.  The curvature of $N$ is always finite (cf., e.g.,
\cite[(2.5)]{Avcurv}), and it is maximal when $N=k$
(\cite[Prop.~2]{Av2}).

The finiteness and extremality statements carry over to the 
``relative" concept defined above, as proved in \cite{AIM}, that is, 
\begin{equation}\label{kextremal}
\cxy_\alpha\!N \leq \cxy_Ak
\qquad\text{and}\qquad
\curv_\alpha\!N \leq \curv_Ak.
\end{equation}

If $N$ happens to be finitely generated over $A$, then its
complexity and curvature are defined over $A$ as well.
In this case relative and absolute notions agree:

\begin{cxycomp} \label{cxycomp}
Let $\alpha\colon A\to B$ be a local homomorphism of local rings. 
If $N$ is finitely generated over $A$, then 
\[
\cxy_\alpha\!N=\cxy_AN \qquad {\textrm{and}}\qquad
\curv_\alpha\!N=\curv_AN.
\]
\end{cxycomp}

This is seen by applying \ref{flcomp} to $M=k$.

\section{The Frobenius endomorphism}\label{frobapp}
\setcounter{subsection}{1}
\setcounter{ssthm}{0}

Let $R$ be a Noetherian ring of characteristic $p$.

We now specialize the results above to the case where $\alpha$ is the
$n$-th Frobenius endomorphism $\varphi^n\colon R\to R$, $M$ is an
$R$-module, and $N=R$.  As explained in Sections \ref{B-structure} and
\ref{A-structure}, we have two $R$-module structures on
$\Tor_i^R(M,\up{\varphi^n}R)$. The structure used is the one that
preserves finiteness properties of $M$, namely the one described in
Section \ref{B-structure}; it comes from the action of $R$ via $\id_R$
on the right-hand variable $R$.  The facts below apply, in particular,
to $F^n(M)=\Tor_0^R(M,\up{\varphi^n}R)$.

The first facts follow from \ref{fg} and \ref{fl}.

\begin{ssprop} \label{R-structureapp} Let $M$ be an $R$-module. 
\begin{enumerate}
\item[\rm{(i)}] 
If $M$ is finitely generated, then so is
$\Tor_i^R(M,\up{\varphi^n}R)$. 
\item[\rm{(ii)}] 
If $M$ has finite length, then so does 
$\Tor_i^R(M,\up{\varphi^n}R)$. 
\qed
\end{enumerate}
\end{ssprop}

Construction \ref{cxcurvextn} yields notions of
complexity and curvature over the powers of the Frobenius endomorphism. 
Inequalities \ref{alphacx}(\ref{kextremal}) yield finiteness statements: 

\begin{ssprop}\label{cxfinitefrob}
Let $R$ be a local ring with residue field $k$.
\[
\cxy_{\varphi^n}\!R\leq \cxy_Rk
\qquad\text{and}\qquad
\curv_{\varphi^n}\!R \leq \curv_Rk.
\]
\end{ssprop}

By \ref{flcomp} and \ref{cxycomp} we have

\begin{ssprop} \label{compfrob}
Let $R$ be as above, and suppose that $\varphi^n$ is finite. 
\begin{enumerate}
\item[\rm{(i)}] 
If $M$ is an $R$-module of finite length, then
\[ 
\ell_R(\Tor_i^R(M,\up{\varphi^n}R))= \ell_{k^{p^n}}(k)\,
\ell_{\varphi^n}(\Tor_i^R(M,\up{\varphi^n}R)). 
\] 
where $k^{p^n}=\varphi^n_k(k)$ is the subfield of $p^n$-th powers.
\item[\rm{(ii)}] There are equalities 
\begin{xxalignat}{3}
&{\phantom{\square}}
&\cxy_{\varphi^n}\!R&=\cxy_R\up{\varphi^n}R
\qquad{\textrm{and}}\qquad
\curv_{\varphi^n}\!R=\curv_R\up{\varphi^n}R\,.
&&\square
\end{xxalignat}
\end{enumerate}
\end{ssprop}

\end{appendix}

\section*{Acknowledgements}

The author wishes to thank Luchezar Avramov whose advice greatly
influenced the paper, both in content and in form, and also Srikanth
Iyengar for interesting and useful comments. In addition, most of the
questions in Section \ref{questions} grew out of conversations
with them.



\begin{thebibliography}{99}

\bibitem{An} 
M.\ Andr\'e, 
\textit{Homologie des Alg\`ebres Commutatives\/}, 
Grundlehren der Math.\ Wiss. {\bf 206}, 
Springer-Verlag, New York-Berlin, 1974.

\bibitem{AB}
M.\ Auslander, M.\ Bridger, 
\textit{Stable Module Theory\/},
Memoirs of the Amer.\ Math.\ Soc.\ {\bf 94} (1969).

\bibitem{Av1}
L.\ L.\ Avramov,
{\em Modules of finite virtual projective dimension\/},
Invent.\ Math.\  {\bf 96} (1989), 71--101.

\bibitem{Avcurv} 
L.\ L.\ Avramov, 
{\em Homological asymptotics of modules over local rings\/}, in 
\textit{Commutative Algebra\/}, 
Math. Sci. Res. Inst. Publ. {\bf 15}
Springer-Verlag, New York-Berlin, 1989, 33--62.

\bibitem{Av2}
L.\ L.\ Avramov,
{\em Modules with extremal resolutions\/},
Math.\ Res.\ Lett.\ {\bf 3} (1996), 319--328.

\bibitem{AFH}
L.\ L.\ Avramov, H.-B.\ Foxby, J.\ Herzog,
{\em Structure of local homomorphisms\/},
J.\ Algebra {\bf 164} (1994), 124--145.

\bibitem{AGP}
L.\ L.\ Avramov, V.\ N.\ Gasharov, I.\ V.\ Peeva,
{\em Complete intersection dimension\/},
Inst.\ Hautes \'Etudes Sci.\ Publ.\ Math.\ {\bf 86} (1997), 67--114.

\bibitem{AIM}
L.\ L.\ Avramov, S.\ Iyengar, C.\ Miller,
{\em Homology of modules over local homomorphisms.
Applications to the Frobenius endomorphism\/},
preprint.

\bibitem{AM}
L.\ L.\ Avramov, C.\ Miller,
{\em Frobenius powers of complete intersections\/},
Math.\ Res.\ Lett.\ {\bf 8} (2001), 225--232.

\bibitem{BFM} P.\ Baum, W.\ Fulton, and R.\ MacPherson, 
{\em Riemann-Roch for singular varieties\/}, 
Inst.\ Hautes \'Etudes Sci.\ Publ.\ Math.\  {\bf 45} (1975), 101--145.

\bibitem{BM}
A.\ Blanco, J.\ Majadas, 
{\em Sur les morphismes d'intersection compl\`ete en caract\'eristique $p$\/},
J.\ Algebra {\bf 208} (1998), 35--42.

\bibitem{BH} 
W.\ Bruns, J.\ Herzog, 
\textit{Cohen-Macaulay Rings\/},
Cambridge Studies in Adv.\ Math.\ {\bf 39},
Cambridge University Press, Cambridge, 1993.

\bibitem{B}
L.\ Burch,
{\em On ideals of finite homological dimension in local rings\/},
Math.\ Proc.\ Cambridge Philos.\ Soc.\  {\bf 64} (1968), 941--952.

 \bibitem{D1}
S.\ P.\ Dutta, 
{\em Frobenius and multiplicities\/},
J.\ Algebra {\bf 85} (1983), 424--448.

\bibitem{D2}
S.\ P.\ Dutta, 
{\em On modules of finite projective dimension over complete intersections\/},
Proc.\ Amer.\ Math.\ Soc.\ {\bf 131} (2003), 113--116. 

\bibitem{DHM} 
S.\ P.\ Dutta, M.\ Hochster, J.\ E.\ McLaughlin, 
{\em Modules of finite projective dimension with negative 
intersection multiplicities\/}, 
Invent.\ Math.\ {\bf 79} (1985), 253--291.

\bibitem{Fu} 
W.\ Fulton, 
\textit{Intersection Theory}, 
Springer-Verlag, New York-Berlin, 1984.

\bibitem{Ge}
A.~A.~Gerko,
{\em On homological dimensions\/}, Mat. Sb. (N.S.) {\bf 192} (2001),
no. 8, 79--94 [Russian]; [English translation: Sb. Math. {\bf 192}
(2001), 1165--1179].

\bibitem{Go}
S.\ Goto, 
{\em A problem on Noetherian local rings of characteristic $p$\/},
Proc.\ Amer.\ Math.\ Soc.\ {\bf 64} (1977), 199--205. 

\bibitem{G}
{T.~H.~Gulliksen},
{\em A change of rings theorem, with applications to Poincar\'e
series and intersection multiplicity\/},
Math. Scand. {\bf 34} (1974), 167--183.

\bibitem{HM} C.\ Han, P.\ Monsky, 
{\em Some surprising Hilbert-Kunz functions\/}, 
Math.\ Z.\ {\bf214} (1993), no. 1, 119--135.

\bibitem{He}
J.\ Herzog, 
{\em Ringe der Charakteristik $p$ und Frobenius-Funktoren\/},
Math.\ Z.\ {\bf 140} (1974), 67--78.

\bibitem{Ho3}
M.\ Hochster,
{\em Cyclic purity versus purity in excellent Noetherian rings\/},
Trans.\ Amer.\ Math.\ Soc.\ {\bf 231} (1977), 464--488.

\bibitem{Hu} 
C.\ Huneke, 
\textit{Tight closure and its applications\/},
Regional Conference Series in Mathematics {\bf 88}, 1996.

\bibitem{HY}
C.\ Huneke, Y.\ Yao,
{\em Unmixed local rings with minimal Hilbert-Kunz multiplicity are regular\/},
Proc.\ Amer.\ Math.\ Soc.\ {\bf 130} (2000), 661--665.

\bibitem{ISW}
S.\ Iyengar, S.\ Sather-Wagstaff, 
{\em The Gorenstein dimension of the Frobenius endomorphism\/},
preprint.

\bibitem{KL1}
J.\ Koh, K.\ Lee,
{\em Some restrictions on the maps in minimal resolutions\/},
J.\ Algebra {\bf 202} (1998), 671--689.

\bibitem{KL2}
J.\ Koh, K.\ Lee,
{\em New invariants of Noetherian local rings\/},
J.\ Algebra {\bf 235} (2001), 431--452.

\bibitem{K}
E.\ Kunz,
{\em Characterization of regular local rings of characteristic $p$\/},
Amer.\ J.\ Math.\ {\bf 41} (1969), 772--784.

\bibitem{Ku} 
K.\ Kurano, 
{\em A remark on the Riemann-Roch formula on affine
schemes associated with Noetherian local rings\/}, 
T\^ohoku Math.\ J.\ {\bf 48} (1996), 121--138.

\bibitem{Kurobring} 
K.\ Kurano, 
{\em On Roberts rings\/}, 
J.\ Math.\ Soc.\ Japan {\bf 53} (2001), 333--355.

\bibitem{Kunumer} 
K.\ Kurano, 
{\em Numerical equivalence defined on a Chow group of a 
Noetherian local ring\/}, 
in preparation.

\bibitem{KS} 
K.\ Kurano, A.\ K.\ Singh, 
{\em Todd classes of affine cones of Grassmannians\/}, 
Int.\ Math.\ Res.\ Notices {\bf 35} (2002), 1841--1855.

\bibitem{L}
C.\ Lech,
{\em Inequalities related to certain couples of local rings\/},
Acta Math.\ {\bf 112} (1964), 69--89.

\bibitem{Mac} 
S.\ Maclane,
\textit{Homology\/},
Grundlehren der Mathematischen Wissenschaften {\bf 114},
Academic Press, New York; Springer-Verlag, New York-Berlin, 1963.

\bibitem{Ma} 
H.\ Matsumura, 
\textit{Commutative Ring Theory\/},
Cambridge Studies in Adv.\ Math.\ {\bf 8},
Cambridge University Press, Cambridge, 1986.

\bibitem{Mi1}
C.\ Miller,
{\em Complexity of tensor products of modules 
and a theorem of Huneke-Wiegand\/},
Proc.\ Amer.\ Math.\ Soc.\ {\bf 126} (1998), 53--60.

\bibitem{Mi}
C.\ Miller,
{\em A Frobenius characterization of finite projective dimension over
complete intersections\/},
Math.\ Z.\ {\bf 233} (2000), 127--136.

\bibitem{MS}
C.\ Miller, A.\ K.\ Singh, 
{\em Intersection multiplicities over Gorenstein rings\/},
Math.\ Ann.\ {\bf 317} (2000), 155--171.

\bibitem{Mo} 
P.\ Monsky, 
{\em The Hilbert-Kunz function\/}, 
Math.\ Ann.\ {\bf 263} (1983), 43--49.

\bibitem{N} 
M.\ Nagata, 
Math.\ Reviews {\bf 40 \# 5609}.

\bibitem{PS1}
C.\ Peskine, L.\ Szpiro,
{\em Sur la topologie des sous-sch\'emas ferm\'es d'un sch\'ema localement 
noeth\'erien, d\'efinis comme support d'un faisceau coh\'erent localement 
de dimension projective finie.\/},
C.\ R.\ Acad.\ Sci.\ Paris S\'er.\ A Math.\  {\bf 269} (1969), 49--51.

\bibitem{PS2}
C.\ Peskine, L.\ Szpiro,
{\em Dimension projective finie et cohomologie locale\/},
Inst.\ Hautes \'Etudes Sci.\ Publ.\ Math.\ {\bf 42} (1973), 47--119.

\bibitem{PS3}
C.\ Peskine, L.\ Szpiro,
{\em Syzygies et multiplicit\'es\/},
C.\ R.\ Acad.\ Sci.\ Paris S\'er.\ A Math.\  {\bf 278} (1974), 1421--1424.

\bibitem{Ro1}
P.\ Roberts, 
{\em The vanishing of intersection multiplicities of perfect complexes\/},
Bull.\ Amer.\ Math.\ Soc.\ (N.S.) {\bf 13} (1985), no.\ 2, 127--130.

\bibitem{RobMacRae}
P.\ Roberts,
{\em The MacRae invariant and the first local Chern character},
Trans.\ Amer.\ Math.\ Soc.\ {\bf 300} (1987), 583--591.

\bibitem{RobMSRI} 
P.\ Roberts, 
{\em Intersection theorems\/}, in 
\textit{Commutative Algebra\/}, 
Math. Sci. Res. Inst. Publ. {\bf 15}
Springer-Verlag, New York-Berlin, 1989, 417--436.

\bibitem{RobBook} 
P.\ Roberts, 
\textit{Multiplicities and Chern Classes in Local Algebra\/}, 
Cambridge Tracts in Mathematics {\bf 133}, 
Cambridge University Press, Cambridge, 1998.

\bibitem{RS}
P.\ Roberts, V.\ Srinivas, 
{\em Modules of finite length and finite projective dimension\/},
Invent.\ Math.\ {\bf 151} (2003), 1--27.

\bibitem{Ro}
A.\ G.\ Rodicio,
{\em On a result of Avramov\/},
Manuscripta Math.\ {\bf 62} (1988), 181--185.

\bibitem{Se}
G.\ Seibert,
{\em Complexes with homology of finite length and Frobenius functors\/},
J.\ Algebra {\bf 125} (1989), 278--287.

\bibitem{Sz}
L.\ Szpiro,
{\em Sur la th\'eorie des complexes parfaits\/}, in
\textit{Commutative algebra\/}, (Durham, 1981), 
London Math.\ Soc.\ Lec.\ Note Ser.\ {\bf 72} (1982), 83--90.

\bibitem{TY}
R.\ Takahashi, Y.\ Yoshino,
{\em Characterizing Cohen-Macaulay local rings by Frobenius maps\/},
preprint.

\bibitem{T}
J.\ Tate,
{\em Homology of Noetherian rings and local rings\/},
Illinois J.\ Math.\ {\bf 1} (1957), 14--27.

\bibitem{WY}
K.\ Watanabe, K.\ Yoshida,
{\em Hilbert-Kunz multiplicity and an inequality between 
multiplicity and colength\/},
J.\ Algebra {\bf 230} (2000), 295--317.

\end{thebibliography}
\end{document}